\documentclass[11pt]{amsart}
\usepackage{amssymb, latexsym, graphicx, mathrsfs, enumerate, amsmath, amsthm, textcomp, url, tikz-cd, nccmath}
\usepackage[T1]{fontenc}
\usepackage{verbatim}
\usepackage[all]{xy}
\usepackage{mathtools}

\usepackage{listings}
\usepackage{dsfont}
\renewcommand{\mathbb}{\mathds}

\usepackage{longtable}

\usepackage[dvipsnames]{xcolor}

\DeclareFontEncoding{OT2}{}{} 
\DeclareTextFontCommand{\textcyr}{\fontencoding{OT2}
    \fontfamily{wncyr}\fontseries{m}\fontshape{n}\selectfont}
\newcommand{\Sha}{{\!\textcyr{Sh}}}

\theoremstyle{plain}
\newtheorem{theorem}{Theorem}[section]

\newtheorem{proposition}[theorem]{Proposition}
\newtheorem{lemma}[theorem]{Lemma}
\newtheorem{corollary}[theorem]{Corollary}

\theoremstyle{definition}

\newtheorem{example}[theorem]{Example}

\newtheorem{remark}[theorem]{Remark}
\newtheorem{definition}[theorem]{Definition}
\newtheorem{construction}[theorem]{Construction}

\newtheorem*{notation*}{Notation}

\usepackage[top=2cm, bottom=3cm, left=3cm, right=3cm]{geometry}
\newcommand{\R}{{\mathbb R}}
\newcommand{\C}{{\mathbb C}}

\newcommand{\Q}{{\mathbb Q}}

\newcommand{\Z}{{\mathbb Z}}

\newcommand{\X}{{\sf X}}

\newcommand{\into}{\hookrightarrow}

\newcommand{\wh}{\widehat}

\DeclareMathOperator{\Aut}{Aut}

\newcommand{\Hom}{{\rm Hom}}

\newcommand{\Gal}{{\rm Gal}}
\newcommand{\ab}{{\rm ab}}

\newcommand{\GL}{{\rm GL}}
\newcommand{\SO}{{\rm SO}}

\newcommand{\ov}{\overline}

\newcommand{\isoto}{\overset{\sim}{\longrightarrow}}
\newcommand{\onto}{\twoheadrightarrow}
\newcommand{\labelto}[1]{\xrightarrow{\makebox[1.5em]{\scriptsize ${#1}$}}}
\newcommand{\labeltoo}[1]{\xrightarrow{\makebox[2.5em]{\scriptsize ${#1}$}}}

\newcommand{\hs}{\kern 0.8pt}
\newcommand{\hssh}{\kern 1.2pt}
\newcommand{\hshs}{\kern 1.6pt}
\newcommand{\hssss}{\kern 2.0pt}

\newcommand{\hm}{\kern -0.8pt}
\newcommand{\hmm}{\kern -1.2pt}

\newcommand{\emm}{\bfseries}

\newcommand{\id}{{\rm id}}

\newcommand{\SL}{{\rm SL}}

\newcommand{\cc}{\raise 1.7pt \hbox{\Tiny{$\bullet$}}}

\newcommand{\coker}{{\rm coker}}

\newcommand{\alg}{{\rm alg}}

\newcommand{\V}{{\mathcal V}}

\newcommand{\Tors}{{\rm Tors}}

\newcommand{\LambdaB}{{\boldsymbol{\Lambda}}}

\newcommand{\Thv}{{\Theta_v}}
\newcommand{\Thvt}{{\Theta_v,\Tors}}
\newcommand{\Tht}{{\Theta,\Tors}}
\newcommand{\ThT}{\Tht}
\newcommand{\vt}{{\vartheta}}

\newcommand{\loc}{{\rm loc}}

\newcommand{\abLambda}{{{\rm ab}\Lambda}}
\newcommand{\mC}{{\mathcal C}}
\newcommand{\lambdatil}{{\tilde\lambda}}

\newcommand{\vn}{{\varnothing}}
\newcommand{\re}{{\rm re}}

\newcommand\Qbar{{\overline\Q}}

\newcommand{\vk}{{\varkappa}}

\newcommand{\la}{\langle}
\newcommand{\ra}{\rangle}
\newcommand{\smm}{{\sf mm}}
\newcommand{\nn}{{\tt nn}}

\newcommand{\Mt}{{\tt M}}
\newcommand{\Mta}{{\tt \la M\ra}}
\newcommand{\Nt}{{\tt N}}
\newcommand{\Nta}{{\tt \la N\ra}}

\newcommand{\other}{{\rm other}}

\newcommand{\vpL}{{\varphi\Lambda}}

\begin{document}

\title[Relative Tate--Shafarevich kernel]%
{Relative Tate--Shafarevich kernel
of a reductive group over  a global field}

\author{Mikhail Borovoi}

\address{School of Mathematical Sciences,
       Tel Aviv University, 6997801 Tel Aviv, Israel}
\email{borovoi@tauex.tau.ac.il}

\thanks{ The author was partially supported
by the Israel Science Foundation (Grant 1030/22).
}



\date{\today}

\begin{abstract}

We describe an algorithm and a program computing the relative Tate--Shafarevich kernel $\Sha^S(K,G)$ 
of a reductive group $G$ over a global field $K$
with respect to a finite set $S$ of places of $K$ containing no archimedean places.
We assume  that we know the algebraic fundamental group $\Lambda=\pi_1^{\rm alg}(G)$
together with the action on it of the Galois group $\Theta=\Gal(L/K)$,
where $L/K$ is a finite Galois extension with no real places
such that $\Gal(K^s/L)$ acts trivially on $\Lambda$.
We also assume that we know the decomposition subgroups
$\Theta_v\subseteq \Theta$ for all places $v$ in $S$,
and all non-cyclic decomposition subgroups $\Theta_v$ for $v$ outside $S$.
Under certain conditions, we compute the set $\Sha^S(K,G)$ also when $S$ contains archimedean places.

\end{abstract}

\maketitle


\section{Introduction}

Let $K$ be a global field, that is, a number field or a global function field, and
let $K^s$ be a fixed separable closure of $K$.
We write  $\V(K)$ or $\V_K$ for the set of places of $K$.
Denote by $\V_\infty(K)$ and $\V_f(K)$ the sets of infinite (archimedean)
and finite (non-archimedean) places of $K$, respectively; then $\V(K)=\V_\infty(K)\cup\V_f(K)$.
Denote by $\V_\R(K)$ and $\V_\C(K)$ the sets of real and complex places of $K$, respectively;
then $\V_\infty(K)=\V_\R(K)\cup \V_\C(K)$.
If $K$ is a function field, then $\V_\infty(K)=\varnothing$ and $\V(K)=\V_f(K)$.

Let $S\subset \V(K)$ be a finite set of places.
Write
\[S_\infty=S\cap\V_\infty(K),\ \, 
S_\R=S\cap\V_\R(K),\ \,  S_\C=S\cap\V_\C(K),\ \, 
S_f=S\cap \V_f(K),\ \, 
S^\complement=\V(K)\smallsetminus S.\]

Let $G$ be a (connected) reductive group over $K$.
We consider  the relative Tate--Shafarevich kernel with respect to $S$:
\[\Sha^S(K,G)\coloneqq\ker \Big[ H^1(K,G)\to\prod_{v\in S^\complement} H^1(K_v,G)\Big], \]
where $H^1(K,G)$ and $H^1(K_v,G)$ denote the corresponding Galois cohomology sets.
This set $\Sha^S(K,G)$ is finite by Borel and Serre \cite[Theorem 7.1]{BS} in the case where $K$ is a number field,
and by Nisnevich \cite{Nisnevich} and by Borel and Prasad \cite[page 176, Theorem]{BP} in the case of a function field.

Note that for $v\in S_\C$ we have $H^1(K_v,G)=1$, and so 
$\Sha^S(K,G)=\Sha^{S\smallsetminus S_\C}(K,G).$
From now on we assume that $S_\C=\vn$; then $S_\infty=S_\R$.

We  refer to \cite[Section 5]{BK} for the definitions of the 
abelian (or abelianized) Galois cohomology groups $H^1_\ab(K,G)$ and $H^1_\ab(K_v,G)$, 
and for the constructions of the canonical abelianization maps
\begin{equation}\label{e:ab-Intro}
\ab\colon H^1(K,G)\to H^1_\ab(K,G),\qquad 
\ab_v\colon H^1(K_v,G)\to H^1_\ab(K_v,G).
\end{equation}
We define an abelian group
\[\Sha^S_\ab(K,G)=\ker\Big[ H^1_\ab(K,G)\to \prod_{v\in S^\complement} H^1_\ab(K_v,G)\Big].\]
Then the abelianization map $\ab$ of \eqref{e:ab-Intro} restricts to an abelianization map
\begin{equation}\label{e:ab-S-Intro}
\ab^S\colon \Sha^S(K,G)\to \Sha^S_\ab(K,G).
\end{equation}

In this paper we show that, under certain assumptions (in
particular, when $S$ contains no archimedean places), the abelianization map \eqref{e:ab-S-Intro} is bijective,
and we compute the finite abelian group $\Sha^S_\ab(K,G)$,
which gives us the cardinality $\#\Sha^S(K,G)=\#\Sha^S_\ab(K,G)$.

Write $\Lambda=\pi_1^\alg(G)$ for the algebraic fundamental group of $G$
of \cite[Section 1]{Borovoi-Memoir}.
Then $\Lambda$ is a finitely generated abelian group, on which
the absolute Galois group $\Gal(K^s/K)$ acts with finite image in $\Aut  \Lambda$.
Let $L$ be a finite Galois extension of $K$ in $K^s$
such that $\Gal(K^s/L)$ acts trivially on $\Lambda$,
and in the case where $K$ is a number field, we choose $L$ to be {\em totally imaginary}.
Then the finite group $\Theta\coloneqq \Gal(L/K)$ naturally acts on $\Lambda$ and on $\V_L$.
(For example, in the case where $K$ is a function field or a totally imaginary number field,
we can take for $L$ the finite Galois extension of $K$
corresponding to the kernel of the canonical homomorphism $\Gal(K^s/K)\to\Aut  \Lambda$;
then $\Theta$ acts on $\Lambda$ effectively, that is, with trivial kernel.)
We shall work with the $\Theta$-module $\Lambda$ and the $\Theta$-set $\V_L$.

For any place $w\in\V_L$, let $\Theta_w\subseteq\Theta$
denote the corresponding decomposition subgroup,
that is, the stabilizer of $w$ in $\Theta$.
Let $v=v(w)\in\V_K$ denote the image of $w$ in $\V_K$.
By abuse of notation, we shall also write $\Theta_v$ for $\Theta_w$.
Then $\Theta_v$ is determined by $v$ up to conjugacy in $\Theta$.

In Section \ref{s:reduction}, using ideas of Tate \cite{Tate} and  of \cite{BK},
we define a finite abelian group $\Sha^S(\Lambda)$ 
in terms of the actions of $\Theta$ on $\Lambda$ and on $\V_L$\hs, as follows.
We consider the $\Theta$-modules
\[\Lambda[\V_L]=\Big\{\sum_{w\in\V_L} \lambda_w\cdot w\  \big|\ \lambda_w\in \Lambda\Big\} \]
(the group of formal finite linear combinations), and
\[ \Lambda[\V_L]_0=\Big\{\sum \lambda_w\cdot w\,\in\, \Lambda[\V_L]\ \big|\ \sum \lambda_w=0\Big\}.\]
We consider the abelian group
$ \big(\Lambda[\V_L]_0\big)_{\Theta,\Tors}$\hs,
the torsion subgroup of the group of coinvariants $\big(\Lambda[\V_L]_0\big)_{\Theta}.$
We define
\begin{equation}\label{e:ShaSLambda}
\Sha^S(\Lambda)=\ker \Big[ \big(\Lambda[\V_L]_0\big)_{\Theta,\Tors}
     \to \bigoplus_{v\in S^\complement} \Lambda_\Thvt\Big].
\end{equation}
See Section \ref{s:reduction} for details.
In that section, we define  canonical maps
\begin{equation}\label{e:ab-Sha-Intro}
\ab\Lambda^S\colon \Sha^S(K,G)\labeltoo{\ \ab^S} \Sha^S_\ab(K,G) \labeltoo{\ \vpL^S} \Sha^S(\Lambda),
\end{equation}
which fit into a canonical commutative diagram
\begin{equation}\label{e:Cart-S-Intro}
\begin{aligned}
\xymatrix@C=14mm{
\Sha^S(K,G)\ar[r]^-{\ab^S} \ar[d]  &\Sha^S_\ab(K,G)\ar[r]^-{\vpL^S} \ar[d]   &\Sha^S(\Lambda)\ar[d] \\
\prod\limits_{v\in S_\infty}\!\! H^1(K_v,G)\ar[r]^-{\ab_{S_\infty}} 
&\bigoplus\limits_{v\in S_\infty}\!\! H^1_\ab(K_v,G)\ar[r]^-{\vpL_{S_\infty}}
     &\bigoplus\limits_{v\in S_\infty}\!\!\Lambda_\Thvt\hskip-4mm
}
\end{aligned}
\end{equation}
The bottom horizontal arrows  
in this diagram 
are the products over $v\in S_\infty$ of the  maps
\begin{equation}\label{e:abLambda-v-Intro} 
\ab\Lambda_v\colon 
H^1(K_v,G)\labeltoo{\ab_v} H^1_\ab(K_v,G)\labeltoo{\vpL_v} \Lambda_\Thvt\hs;
\end{equation}
see Section \ref{s:reduction}. We prove the following theorem.

\begin{theorem}[Theorem \ref{t:Sha}]
\label{t:Sha-Intro}
Both rectangles in diagram \eqref{e:Cart-S-Intro} are Cartesian. 
\end{theorem}

We deduce two corollaries:

\begin{corollary}[Corollary \ref{c:bijective}]
\label{c:bijective-Intro}
If for all $v\in S_\infty=S_\R$  the   maps $\ab\Lambda_v$
of \eqref{e:abLambda-v-Intro} are bijective,
then 
both arrows in \eqref{e:ab-Sha-Intro} are bijective.
\end{corollary}

\begin{corollary}[Corollary \ref{c:no-infinite}]
\label{c:no-infinite-Intro}
If $S$ contains no archimedean  places, then
both arrows in \eqref{e:ab-Sha-Intro} are bijective.
In particular, this is the case when $K$ is a global function field.
\end{corollary}

Note that when $S$ is as in Corollary \ref{c:bijective-Intro}, 
the bijection $\ab^S$ of \eqref{e:ab-Sha-Intro} 
induces on the set $\Sha^S(K,G)$ a canonical abelian group structure.

The abelian group $\Sha^S(\Lambda)$ is defined in \eqref{e:ShaSLambda}
in terms of {\em infinite} sets $\V_L$ and $S^\complement$.
For a computer computation, we must describe $\Sha^S(\Lambda)$ in terms of finite sets.

Let $\Sigma$ be a finite or infinite set such that  $S\subseteq\Sigma\subseteq\V_K$.
We write $\Sigma_L$ for the preimage of $\Sigma$ in $\V_L$, 
and we consider the $\Theta$-module
\[ \Lambda[\Sigma_L]_0=\Big\{\sum \lambda_w\cdot w\,\in\, \Lambda[\Sigma_L]\ \big|\ \sum \lambda_w=0\Big\}.\]
We consider the abelian group
$ \big(\Lambda[\Sigma_L]_0\big)_{\Theta,\Tors}$\hs,
the torsion subgroup of the group of coinvariants $\big(\Lambda[\Sigma_L]_0\big)_{\Theta}.$
We define
\begin{equation*}
\ker_\Sigma^S(\Lambda)=\ker \Big[\big(\Lambda[\Sigma_L]_0\big)_\Tht\to \bigoplus_{v\in \Sigma\smallsetminus S}\Lambda_\Thvt\Big].
\end{equation*}
In particular, if $\Sigma=\V_K$, then $\ker_\Sigma^S(\Lambda)=\Sha^S(\Lambda)$.

\begin{definition}
We say that a set  $\Sigma$ of places of $K$ satisfying $S\subset\Sigma$ is {\em sufficiently large for $S$}
if for any $v\in \Sigma^\complement$ there exists $v_\Sigma\in \Sigma\smallsetminus S$
such that $\Theta_v\subseteq\Theta_{v_\Sigma}$ (up to conjugacy).
\end{definition}

By Lemma \ref{l:S-Sigma} below, for any finite set $S\subset\V_K$ there exists a {\em finite} set $\Sigma\supset S$ 
that is sufficiently large for $S$. Note that when $S$ is finite and $\Sigma$ sufficiently large for $S$,
the inclusion $S\subset \Sigma$ is proper; this is why we write $S\subset \Sigma$ rather than $S\subseteq \Sigma$.

\begin{theorem}[Theorem \ref{t:Sigma'}] 
\label{t:Sigma'-Intro}
Let $S\subset\Sigma\subseteq\Sigma'\subseteq\V_K$
be subsets, where $\Sigma$ is sufficiently large for $S$.
Then the natural homomorphism
\[\ker^S_\Sigma(\Lambda)\to\ker^S_{\Sigma'}(\Lambda)\]
is an isomorphism.
In particular  (taking $\Sigma'=\V_K$), the natural homomorphism $\ker^S_\Sigma(\Lambda)\to\Sha^S(\Lambda)$ is an isomorphism.
\end{theorem}

Theorem \ref{t:Sigma'-Intro} together with Corollaries \ref{c:bijective-Intro} and  \ref{c:no-infinite-Intro}
permits us to compute (under certain conditions) the finite abelian group $\Sha^S_\ab(K,G)$  using a computer.
Note that Theorems \ref{t:Sha-Intro} and \ref{t:Sigma'-Intro} together 
give an alternative proof of the finiteness of $\Sha^S(K,G)$
(we use the finiteness of $H^1(K_v,G)$ for any real place $v$ of $K$; see, 
for instance, \cite{Borovoi-real}\hs).

Note that there are computations of  $\Sha(K,G)\coloneqq\Sha^\varnothing(K,G)$
for certain connected reductive $K$-groups in \cite[Examples 5.6, 5.7, 5.8]{Sansuc}, \hs
\cite[\,(4.2.2)\hs]{Kottwitz-84}, \hs\cite[Corollary 5.13]{Borovoi-Memoir}, \cite[Theorem 9.4]{CT}, and
\cite[Theorem 1.10]{Rosengarten:Tate}.
Moreover, Sansuc \cite{Sansuc} computes $\Sha^S(K,B)$ with {\em non-empty} $S$
for some {\em finite} abelian $K$-groups $B$;
see \cite[Sections 1,\,2]{Sansuc}, where he also computes $\Sha_\omega(K,B)\coloneqq\varinjlim_S \Sha^S(K,B)$.
However, no computations of $\Sha^S(K,G)$ for connected reductive $K$-groups $G$ and non-empty $S$
in the literature are known to the author (except for the classical cases where $G$ is simply connected semisimple, of $G=\GL_n$, 
or $G=\SO(K^n,{\mathcal F})$ for a quadratic form $\mathcal F$ in $n$ variables).

The plan of the rest of the paper is as follows.
In Section \ref{s:reduction}  we prove Theorem \ref{t:Sha-Intro}, from which Corollaries 
\ref{c:bijective-Intro} and \ref{c:no-infinite-Intro} follow.
These corollaries reduce (under certain conditions)  computing the set $\Sha^S(K,G)$ to computing the abelian group $\Sha^S_\ab(K,G)\cong\Sha^S(\Lambda)$.
In Section \ref{s:abelian} we prove Theorem \ref{t:Sigma'-Intro},
which describes $\Sha^S(\Lambda)$ in terms of finitely many operations 
with finitely generated abelian groups.
In Section \ref{s:algorithm} we describe our algorithm and the Magma program,
which are based on the results of Sections \ref{s:reduction} and \ref{s:abelian}.
In Section \ref{s:examples} we consider a certain semisimple group $G$ over $\Q$
that splits over $\Q(\zeta_8)$ where $\zeta_8$ is a primitive 8th root of unity.
We give  examples of computations of the finite abelian groups
$\Sha^S(\Lambda)\cong\Sha^S_\ab(K,G_K)$ 
for the fields $K=\Q\big(\sqrt7\big)$ and $K=\Q(\alpha)$ where $\alpha$
is a root of the polynomial $f(x)=x^3+x^2-10x-8$. 
In Appendix \ref{app:simple} we describe  a simple alternative method 
(suggested by Zev Rosengarten) of computing the cardinality of the finite set $\Sha^S(K,G)$ 
for $S$ as in Corollary \ref{c:bijective-Intro}.
In general, this alternative method does not give 
the isomorphism class of the finite abelian group $\Sha^S_\ab(K,G)$.
In this arXiv version of the paper, in  Appendix B
we provide a listing of our Magma program.

{\sc Acknowledgments\ } We thank Willem A. de Graaf 
for writing the Magma function $\tt ComplementIntMat$ used in our Magma program. 
We thank David Harari, Taeyeoup Kang, and Zev Rosengarten 
for helpful discussions and email correspondence.
We used the AI assistant Claude (Anthropic) to proofread the paper and to check the program.

\section{Reduction to computing $\Sha^S(\Lambda)$}
\label{s:reduction}

In this section we reduce computing $\Sha^S(K,G)$
to computing a certain abelian group $\Sha^S(\Lambda)$.

Recall that for $v\in \V(K)$, we denote by
$\Thv\subseteq \Theta$  a decomposition subgroup of $\Theta$ at $v$
(determined up to conjugacy).
We write $\Lambda_\Thv$ for the group of coinvariants of $\Thv$ in $\Lambda$, that is,
\[\Lambda_\Thv=\Lambda/\langle \vartheta\cdot x-x\ | \ x\in \Lambda,\, \vartheta\in\Thv\rangle.\]
We write
\[\Lambda_{\Thv,\Tors}\coloneqq (\Lambda_\Thv)_\Tors\]
for the torsion subgroup of $\Lambda_\Thv$.

Let $\Sigma\subseteq\V_K$ be a subset, finite or infinite.
We write $\Sigma_L$ for the preimage of $\Sigma$ in $\V_L$.
We consider the $\Theta$-modules
\[\Lambda[\Sigma_L]=\Big\{\sum_{w\in\Sigma_L} \lambda_w\cdot w\  \big|\ \lambda_w\in \Lambda\Big\} \]
(the group of formal finite linear combinations), and
\[ \Lambda[\Sigma_L]_0=\Big\{\sum \lambda_w\cdot w\,\in\, \Lambda[\Sigma_L]\ \big|\ \sum \lambda_w=0\Big\}.\]
The Galois group $\Theta$ acts on $\Lambda[\Sigma_L]$ and $\Lambda[\Sigma_L]_0$
diagonally: for $\vt\in\Theta$ we have
\[\vt\cdot \sum_{w\in\Sigma_L}\! \lambda_w\cdot w\hs=
\sum_{w\in\Sigma_L}\!(\vt\cdot \lambda_w)\cdot (\vt\cdot w)
\hs=\sum_{w\in\Sigma_L}\!\!\big(\vt\cdot \lambda_{\vt^{-1}\cdot w}\big)\cdot w.
\]
We consider the abelian group
\[ \big(\Lambda[\Sigma_L]_0\big)_{\Theta,\Tors}\hs,\]
the torsion subgroup of the group of coinvariants $\big(\Lambda[\Sigma_L]_0\big)_{\Theta}.$

For a place $v\in\Sigma$, let $\V_L(v)$ denote the preimage of $v$ in $\V_L$ (hence in $\Sigma_L$).
We have a canonical homomorphism
\[\Lambda[\Sigma_L]\to \Lambda[\V_L(v)], \quad\
   \sum_{w\in \Sigma_L}\lambda_w\cdot w\,\mapsto\, \sum_{w\in \V_L(v)} \lambda_w\cdot w.\]
Choose $w\in\V_L(v)$.
Since the $\Theta$-module $\Lambda[\V_L(v)]$ is induced from the $\Theta_w$-module $\Lambda$, by Shapiro's lemma
(see \cite[Proposition 9.76(ii)]{Rotman})
the embedding
\[ \Lambda\into \Lambda[\V_L(v)],\quad\ \lambda\mapsto \lambda\cdot w \]
induces an isomorphism
\[\Lambda_{\Theta_v}= H_0(\Theta_v,\Lambda)\isoto H_0(\Theta,\Lambda[\V_L(v)])= \Lambda[\V_L(v)]_\Theta\hs,\]
whence we obtain an isomorphism
\[ \Lambda_\Thvt\isoto \Lambda[\V_L(v)]_\Tht\hs.\]
Thus we obtain a canonical homomorphism
\begin{equation}\label{e:lv}
l_v\colon\,\big(\Lambda[\Sigma_L]_0\big)_{\Theta,\Tors}
    \to \Lambda[\Sigma_L]_{\Theta,\Tors}\to
   \Lambda[\V_L(v)]_{\Theta,\Tors}\overset\sim\longrightarrow \Lambda_{\Thv,\Tors}
\end{equation}
where the last arrow is inverse to the isomorphism above.

We refer to \cite[Chapter 5]{BK}
for the definition of the abelian Galois cohomology groups
$H^1_\ab(K,G)$ and $H^1_\ab(K_v,G)$
and the construction of the abelianization maps
\[ \ab\colon H^1(K,G)\to H^1_\ab(K,G)\quad\ \text{and}
    \quad\ \ab_v\colon H^1(K_v,G)\to H^1_\ab(K_v,G).\]

For a {\em real} place $v$ of $K$, 
the abelianization map
\[\ab_v\colon H^1(K_v,G)\to H^1_\ab(K_v,G)\]
is surjective; see \cite[Theorem 5.4]{Borovoi-Memoir}.
Since $L$ is totally imaginary, we have
$L_w\simeq\C$ for every place $w$ of $L$ over $v$, hence $\#\Theta_v=2$.
We have the Tate--Nakayama isomorphism
\begin{equation*}
H^1_\ab(K_v,G)\cong H^{-1}(\Theta_v,\Lambda)
\end{equation*}
where $H^{-1}$ denotes the corresponding Tate cohomology group;
see \cite[Proposition 8.21]{BT}.
By definition, we have
\[H^{-1}(\Theta_v,\Lambda)=Z^{-1}(\Thv,\Lambda)/B^{-1}(\Thv,\Lambda)\]
where 
\[Z^{-1}(\Thv,\Lambda)=\{x\in \Lambda \ |\ x+\vt\cdot x=0\}
\quad\text{and}\quad B^{-1}(\Thv,\Lambda)= \{y-\vt\cdot y\ |\ y\in \Lambda\},\]
$\vt$ denoting the non-trivial element of $\Theta_v$.
Since $\Lambda/B^{-1}(\Thv,\Lambda)=\Lambda_\Thv$, we see that 
$H^{-1}(\Theta_v,\Lambda)\subseteq \Lambda_\Thv$. 
For any $x\in Z^{-1}(\Thv,\Lambda)$ we have
\[[2x] =[x+\vt\cdot x] + [x-\vt\cdot x]=0+0=0\in H^{-1}(\Theta_v,\Lambda),\]
where $x+\vt\cdot x=0$ because $x\in Z^{-1}(\Thv,\Lambda)$,
and where
$[x-\vt\cdot x]=0\in H^{-1}(\Thv,\Lambda)$ 
because $x-\vt\cdot x\in B^{-1}(\Thv,\Lambda)$.
It follows that $H^{-1}(\Theta_v,\Lambda)\subseteq \Lambda_\Thvt$.
Thus we obtain an injective homomorphism
\begin{equation*}
\vpL_v\colon H^1_\ab(K_v,G)\cong H^{-1}(\Theta_v,\Lambda)\into \Lambda_\Thvt\hs.
\end{equation*}
We obtain a map
\begin{equation}\label{e:abLambda-v}
\ab\Lambda_v\colon H^1(K_v,G)\labelto{\ab_v}  H^1_\ab(K_v,G)
\labelto{\vpL_v}\Lambda_{\Theta_v,\Tors}
\end{equation}
where the map $\ab_v$ is surjective, and the homomorphism $\vpL_v$
is injective. It follows that if the composite map $\abLambda_v$
is bijective, then so are both $\ab_v$ and $\vpL_v$.

For a {\em complex} place $v$  of $K$ we have the trivial embedding
\[\ab\Lambda_v\colon H^1(K_v,G)=H^1_\ab(K_v,G)=\{1\}\into \Lambda_\Thvt=\Lambda_\Tors\hs.\]

For a {\em finite} place $v$ of $K$, we have  canonical bijections
\[\ab\Lambda_v\colon H^1(K_v,G)\labelto{\ab_v} H^1_\ab(K_v,G)
\labelto{\vpL_v} \Lambda_{\Theta_v,\Tors};\]
see \cite[Theorem 5.5.2(1)]{BK}.

For our global field $K$, we have  canonical maps
\begin{equation}\label{e:ab}
\abLambda\colon H^1(K,G)\labelto\ab H^1_\ab(K,G)\labelto{\vpL} \big(\Lambda[\V_L]_0\big)_\Tht\hs;
\end{equation}
for the construction of the homomorphism $\vpL$ see \cite[Theorem 5.5.4(1) and Definition 4.3.7]{BK}.

Set
\[\LambdaB_{\Sigma_\infty}=\bigoplus_{v\in\Sigma_\infty} \Lambda_{\Thv,\Tors}\]
where $\Sigma_\infty=\Sigma\cap\V_\infty(K)$.
We have a canonical homomorphism
\[\bigoplus_{v\in\Sigma_\infty} l_v \colon\,\big(\Lambda[\Sigma_L]_0\big)_{\Theta,\Tors}\to
      \bigoplus_{v\in\Sigma_\infty} \Lambda_{\Thv,\Tors}=\LambdaB_{\Sigma_\infty}\hs.\]
On the other hand, using the maps $\ab\Lambda_v$ for $v\in\Sigma_\infty$ defined above,
we obtain a canonical map
\[ \ab\Lambda_{\Sigma_\infty}\colon\prod_{v\in \Sigma_\infty} H^1(K_v,G)\to \LambdaB_{\Sigma_\infty}\hs.\]

Let $S$ be a subset of $\V_K$ and write $S_L$ for the preimage of $S$ in $\V_L$\hs.
Write $S^\complement=\V_K\smallsetminus S$ for the complement of $S$ in $\V_K$\hs,
and $S_L^\complement=\V_L\smallsetminus S_L$ for the complement of $S_L$ in $\V_L$.
Recall that 
\begin{align*}
\Sha^S(\Lambda)  &= \ker \Big[\bigoplus_{v\in S^\complement}l_v\colon\hs  \big(\Lambda[\V_L]_0\big)_{\Theta,\Tors}\to \bigoplus_{v\in S^\complement} \Lambda_\Thvt\,\Big]\\
&=\ker \Big[ \big(\Lambda[\V_L]_0\big)_{\Theta,\Tors}
     \to \Lambda\big[S_L^\complement\big]_{\Theta,\Tors}\Big],
\end{align*}
where
\[\Lambda\big[S_L^\complement\big]_{\Theta,\Tors}
   \cong \hs \bigoplus_{v\in S^\complement} \Lambda_\Thvt\quad\  \text{because} \quad\ \Lambda\big[S_L^\complement\big]=\bigoplus_{v\in S^\complement}\Lambda[\V_L(v)].\]
The maps \eqref{e:ab} induce  maps
\begin{equation}\label{e:ab-Sha}
\abLambda^S\colon\, \Sha^S(K,G)\labelto{\ab^S} \Sha^S_\ab(K,G)
\labelto{\vpL^S} \Sha^S(\Lambda)
\end{equation}
(see the proof of Theorem \ref{t:Sha} below).
Moreover, we have a natural map
\[ l_{S_\infty}=\bigoplus_{v\in S_\infty}l_v\colon\hs \Sha^S(\Lambda)\to \bigoplus_{v\in S_\infty} \Lambda_\Thvt\eqqcolon\LambdaB_{S_\infty}\hs.\]

\begin{theorem}\label{t:Sha}
Let $S$ be a subset of $\V_K$\hs, finite or infinite.
Both  rectangles in the following diagram  are  Cartesian:
\begin{equation}\label{e:Cart-S}
\begin{aligned}
\xymatrix@C=14mm{
\Sha^S(K,G)\ar[r]^-{\ab^S} \ar[d]^-{\loc_{S_\infty}}
&\Sha^S_\ab(K,G)\ar[r]^-{\vpL^S} \ar[d]^-{\loc_{S_\infty}}
&\Sha^S(\Lambda)\ar[d]^-{l_{S_\infty}} \\
\prod\limits_{v\in S_\infty}\!\!\! H^1(K_v,G)
\ar[r]^-{\ab_{S_\infty}}  
&\bigoplus\limits_{v\in S_\infty}\!\!\! H^1_\ab(K_v,G)
\ar[r]^-{\vpL_{S_\infty}}  
&\LambdaB_{S_\infty}\hskip-5mm
}
\end{aligned}
\end{equation}
\end{theorem}

\begin{proof}
Diagram \eqref{e:Cart-S} naturally embeds into the following diagram:
\begin{equation}\label{e:Cart}
\begin{aligned}
\xymatrix@C=14mm{
H^1(K,G)\ar[r]^-{\ab} \ar[d]^-{\loc_\infty}
&H^1_\ab(K,G)\ar[r]^-{\vpL} \ar[d]^-{\loc_\infty}
&\big(\Lambda[\V_L]_0\big)_\Tht\hskip-7mm\ar[d]^-{l_\infty} \\
\prod\limits_{v\in \V_\infty(K)}\!\! H^1(K_v,G)\ar[r]^-{\ab_{\V_\infty(K)}}
&\bigoplus\limits_{v\in \V_\infty(K)}\!\! H^1_\ab(K_v,G)\ar[r]^-{\vpL_{\V_\infty(K)}}
       &\LambdaB_{\V_\infty(K)}\hskip-8.5mm
}
\end{aligned}
\end{equation}
in which all rectangles are Cartesian by the proof of \cite[Theorem 5.5.10]{BK}.
It follows that diagram \eqref{e:Cart-S} is commutative.

We prove that the right-hand rectangle in diagram \eqref{e:Cart-S} is Cartesian.
We must show that if elements $x_{S_\infty}\in \bigoplus_{v\in S_\infty}\! H^1_\ab(K_v,G)$ and
$x_\Lambda\in \Sha^S(\Lambda)$ have the same image in
$\LambdaB_{S_\infty}$,
then they are the images of a unique element
$x\in\Sha^S_\ab(K,G)$.

Let
\[ x_\infty\in  \bigoplus\limits_{v\in \V_\infty(K)}\!\! H^1_\ab(K_v,G)\]
be the extension of $x_{S_\infty}$ by the trivial classes
at $v\in\V_\infty(K)\smallsetminus S_\infty$.
Then the elements $x_{\infty}=(x_v)_{v\in \V_\infty(K)}$ and $x_\Lambda$ have the same image in $\LambdaB_{\V_\infty(K)}$, that is, the image of $x_v$
in $\Lambda_\Thvt$ is equal to $l_v(x_\Lambda)$.
Indeed, for $v\in S_\infty$ this is the hypothesis, and 
for each archimedean $v\notin S$ we have $l_v(x_\Lambda)=0$,
which holds because $x_\Lambda\in\Sha^S(\Lambda)$.
Since the right-hand rectangle in diagram \eqref{e:Cart} is Cartesian,
the elements $x_{\infty}$ and $x_\Lambda$ are the images of a unique element $x\in H^1_\ab(K,G)$,
and it remains to show that $x\in\Sha^S_\ab(K,G)$.
By construction, for each {\em infinite} $v\in \V_\infty(K)\smallsetminus S_\infty$\hs,
the image of $x$ in  $H^1_\ab(K_v,G)$ is trivial.
Moreover, for each {\em finite} place $v\in\V_f(K)\smallsetminus S_f$\hs,
by \cite[the proof of Proposition 4.6.2]{BK}
the diagram below is commutative and
the bottom horizontal arrow in this diagram is an isomorphism:
\[
\xymatrix@C=14mm{
H^1_\ab(K,G)\ar[r]^-\vpL\ar[d]_-{\loc_v} & \big(\Lambda[\V_L]_0\big)_\Tht\hskip-7mm\ar[d]^-{l_v}\\
H^1_\ab(K_v,G) \ar[r]^-{\vpL_v}  &\Lambda_\Thvt\,.\hskip-10mm
}
\]
Therefore, for $l_v$ as in \eqref{e:lv}, 
from $l_v(\vpL(x))=0$ (because $x_\Lambda\in\Sha^S(\Lambda)$\hs)  we conclude that $\loc_v(x)=1$.
Thus for all
\[v\in(\V_\infty(K)\smallsetminus S_\infty)\cup(\V_f(K)\smallsetminus S_f)=S^\complement\]
we have $\loc_v(x)=1$, that is, $x\in\Sha^S_\ab(K,G)$.
We see that  the right-hand rectangle in diagram \eqref{e:Cart-S} is Cartesian. 

We prove similarly that the left-hand rectangle in diagram \eqref{e:Cart-S} is Cartesian as well, 
using the fact that for each   {\em finite} place $v\in\V_f(K)\smallsetminus S_f$\hs,
the diagram below is commutative (see above) and
the bottom horizontal arrow in this diagram is bijective:
\[
\xymatrix@C=14mm{
H^1(K,G)\ar[r]^-\ab\ar[d]_-{\loc_v} &H^1_\ab(K,G)\ar[d]^-{\loc_v}\\
H^1(K_v,G) \ar[r]^-{\ab_v}  &H^1_\ab(K_v,G)\hs,
}
\]
which completes the proof of the theorem.
\end{proof}

\begin{corollary}[Corollary \ref{c:bijective-Intro}]
\label{c:bijective}
If for all $v\in S_\infty$  the  maps $\ab\Lambda_v$
of \eqref{e:abLambda-v} are bijective,
then both arrows in  \eqref{e:ab-Sha} are bijective.
\end{corollary}

\begin{proof} Indeed, then both bottom horizontal arrows 
in diagram \eqref{e:Cart-S} with Cartesian rectangles are bijective
(see the discussion after \eqref{e:abLambda-v}\hs).
Since a base change of a bijection is a bijection,
both top horizontal arrows are bijective as well.
\end{proof}

\begin{corollary}[Corollary \ref{c:no-infinite-Intro}]
\label{c:no-infinite}
If $S$ contains no infinite (archimedean)  places,
then both arrows in  \eqref{e:ab-Sha} are bijective.
In particular, this is the case when $K$ is a global function field.
\end{corollary}

\begin{proof}
Indeed, then $S_\infty=\varnothing$, and the hypothesis of Corollary \ref{c:bijective}
is vacuous.
\end{proof}

\section{$\Sha^S(\Lambda)$ in terms of finite abelian groups}
\label{s:abelian}

By Theorem \ref{t:Sha}, in order to compute $\Sha^S(K,G)$ for a finite subset $S\subset \V(K)$, 
it suffices to compute the finite abelian group $\Sha^S(\Lambda)$ and, 
for each {\em real} place $v\in S_\infty$, the finite pointed set
$H^1(K_v,G)$ together with the maps
\[\abLambda_v\colon\,  H^1(K_v,G)\to H^1_\ab(K_v,G)\cong H^{-1}(\Theta_v,\Lambda)\into  \Lambda_\Thvt\]
and
\[ l_v\colon \Sha^S(\Lambda)\to \Lambda_\Thvt\hs.\]
The real cohomology sets $H^1(K_v,G)$ can be computed using the computer program of \cite{BdG}.
This program gives explicit cocycles in $Z^1(K_v,G)$,
and so it is not difficult to compute their images in
$\Lambda_\Thvt$ under $\abLambda_v$; we shall describe this in detail elsewhere.

It remains to compute the finite abelian group $\Sha^S(\Lambda)$
together with the homomorphisms $l_v$ for $v\in S_\infty$.
In this section we give a formula for $\Sha^S(\Lambda)$
in terms of certain  finite abelian groups, and in Section \ref{s:algorithm}
we give formulas for $l_v$.
This permits one to compute $\Sha^S(\Lambda)$ and the homomorphisms $l_v$ using a computer.

For a finite Galois extension of global fields $L/K$ with Galois group $\Theta=\Gal(L/K)$,
and  a finite or infinite  set of places $S\subset \V_K$,
we consider a finite or infinite set $\Sigma\subseteq\V(K)$ containing $S$.
We denote by $S_L$ and $\Sigma_L$ the preimages of  $S$ and $\Sigma$ in $\V(L)$,
respectively.

\begin{definition}\label{d:sl}
For $S\subset\Sigma\subseteq \V(K)$,
we say that  $\Sigma$ is {\em sufficiently large for $S$}
if for any place $w\in \Sigma_L^\complement$
there exists a  place $w_{_\Sigma}\in \Sigma_L\smallsetminus S_L$
such that
\begin{equation*}
\Theta_{w_{_\Sigma}}\supseteq \Theta_w\hs.
\end{equation*}
\end{definition}

\begin{lemma} \label{l:S-Sigma}
For any finite set $S\subset \V_K$,
there exists a {\emm finite} set $\Sigma\subset \V(K)$ containing $S$
that is sufficiently large for $S$.
\end{lemma}

\begin{proof}
Let $\V_{\rm nc}(L/K)\subset\V(K)$ denote the set of places of $K$
with non-cyclic decomposition group in $L$.
At each place $v\in\V_{\rm nc}(L/K)$, the finite extension $L/K$ ramifies,
and therefore there are only finitely many such places.
Thus we may choose $\Sigma$ containing the finite sets $S$ and  $\V_{\rm nc}(L/K)$.
Since our Galois group $\Theta=\Gal(L/K)$ is finite,
there are finitely many conjugacy classes $\mC$ of maximal cyclic subgroups $C\subseteq\Theta$.
For each such conjugacy class $\mC$ we choose $C\in\mC$.
Since $S_L$ is finite, by the Chebotarev density theorem
we can choose a finite place $w_{C}\in (\V_L\smallsetminus S_L)_f$ with decomposition group $C$.
Let $v_\mC \in (\V_K\smallsetminus S)_f$ denote the image of $w_C$.
We may choose $\Sigma$ containing the finite set of these $v_\mC$.

Now for $\Sigma$ as above and for any $w\in \Sigma_L^\complement$,
the decomposition group $\Theta_w$ of $w$ is cyclic,
and hence $\Theta_w\subseteq C_0$ for some maximal cyclic subgroup $C_0$ of $\Theta$.
By construction, there exists a finite place $w_C\in \Sigma_L\smallsetminus S_L$
with decomposition group $\Theta_{w_C}=\vt C_0\vt^{-1}$ for some $\vt\in\Theta$.
Set $w'= \vt^{-1}\cdot w_C$;
then $\Theta_{w'} =C_0$.
Since both $\Sigma_L$ and $S_L$ are $\Theta$-invariant, we have
$w'\in\Sigma_L\smallsetminus S_L$.
We obtain that $\Theta_w\subseteq C_0=\Theta_{w'}$ for $w'\in \Sigma_L\smallsetminus S_L$\hs.
Thus  $\Sigma$ is sufficiently large for $S$ and finite, as desired.
\end{proof}

\begin{construction}
Let $S\subset \V_K$ be a  set of places of $K$,
and let $\Sigma\subseteq\V_K$ be a sufficiently large set for $S$, finite or infinite.
We construct a $\Theta$-equivariant map
\begin{equation*}
\phi_\Sigma\colon \Sigma_L^\complement \to\Sigma_L\smallsetminus S_L
\end{equation*}
as follows.
For any $v\in \Sigma^\complement$
we choose a place $w\in \Sigma_L^\complement$ over $v$.
Then there exists a place $w_\Sigma\in \Sigma_L\smallsetminus S_L$
with $\Theta_{w_\Sigma}\supseteq \Theta_w$.
We define the map $\phi_\Sigma$ on the orbit $\V_L(v)=\Theta\cdot w$
by $\phi_\Sigma(\vartheta\cdot w)=\vartheta\cdot w_\Sigma$ for $\vartheta\in\Theta$.
Since $\Theta_{w_\Sigma}\supseteq \Theta_w$\hs,
this map is well defined.

We extend the map $\phi_\Sigma$ to a map
\begin{equation}\label{e:phi-ext}
\wh\phi_\Sigma\colon \V_L \to \Sigma_L
\end{equation}
by putting $\wh\phi_\Sigma(w)=w$ for $w\in \Sigma_L$. Let
\[\iota\colon \Sigma_L\into\V_L\]
denote the inclusion map.
Then
\begin{equation*}
\wh\phi_\Sigma\circ\iota=\id.
\end{equation*}
By abuse of notation, we write $\phi_\Sigma$ for $\wh\phi_\Sigma$.
Then the map $\phi_\Sigma\colon \V_L\to\Sigma_L$ is surjective.
\end{construction}

Let $\Lambda$ be a $\Theta$-module.
For a finite or infinite set $\Sigma$ satisfying $S\subseteq\Sigma\subseteq\V_K$,
we define
\begin{multline*}
\ker_\Sigma^S(\Lambda)=\ker \Big[\big(\Lambda[\Sigma_L]_0\big)_\Tht\to \bigoplus_{v\in \Sigma\smallsetminus S}\Lambda_\Thvt\Big]\\
   =\ker \Big[\big(\Lambda[\Sigma_L]_0\big)_\Tht\to \Lambda[\Sigma_L\smallsetminus S_L]_\Tht\Big]
\end{multline*}
where the homomorphism $\big(\Lambda[\Sigma_L]_0\big)_\Tht\to \Lambda[\Sigma_L\smallsetminus S_L]_\Tht$ is given by restriction from $\Sigma_L$ to $\Sigma_L\smallsetminus S_L$\hs:
\[ \bigg[\,\sum_{w\in \Sigma_L} \lambda_w\cdot w\hs\bigg]\longmapsto
    \bigg[\,\sum_{w\in \Sigma_L\smallsetminus S_L}\!\!\!\! \lambda_w\cdot w\hs\bigg].\]
Then
\[ \ker^S_{\V(K)}(\Lambda)=\Sha^S(\Lambda).\]

\begin{theorem}[Theorem \ref{t:Sigma'-Intro}]
\label{t:Sigma'}
Let $L/K$ be a finite Galois extension of global fields with Galois group $\Theta=\Gal(L/K)$,
and let $\Lambda$ be a $\Theta$-module.
Let
\[ S\subset\Sigma\subseteq\Sigma'\subseteq\V_K\]
be subsets, where $\Sigma$ is sufficiently large for $S$.
Then the natural homomorphism
\[\ker^S_\Sigma(\Lambda)\to\ker^S_{\Sigma'}(\Lambda)\]
is an isomorphism. 
In particular (taking $\Sigma'=\V_K$), the natural homomorphism
$\ker^S_\Sigma(\Lambda)\to\Sha^S(\Lambda)$ is an isomorphism.
\end{theorem}

\begin{proof}[Proof of Theorem \ref{t:Sigma'}]
We denote by $\phi'_\Sigma\colon \Sigma'_L\to\Sigma_L$ the restriction to $\Sigma'_L$
of the map $\phi_\Sigma \colon \V_L\to \Sigma_L$ of \eqref{e:phi-ext}.
To simplify the notation, we write
$\phi$ for $\phi'_\Sigma$.
From our construction of the map $\phi$, it follows that
$\phi(\Sigma'_L\smallsetminus S_L)\subseteq \Sigma_L\smallsetminus S_L$.
From the following commutative diagram:
\[
\xymatrix{
\Lambda[\Sigma'_L]_0\ar[r]\ar[d]^{\phi_*} &\Lambda[\Sigma'_L\smallsetminus S_L]\ar[d]^{\phi_*}\\
\Lambda[\Sigma_L]_0\ar[r]                &\Lambda[\Sigma_L\smallsetminus S_L]
}
\]
we obtain a commutative diagram with exact rows
\[
\xymatrix{
0\ar[r] &\ker^S_{\Sigma'}\ar@{-->}[d]^{\phi_*}\ar[r]
        &\big(\Lambda[\Sigma'_L]_0\big)_\Tht\ar[r]\ar[d]^{\phi_*}
       &\Lambda[\Sigma'_L\smallsetminus S_L]_\Tht\ar[d]^{\phi_*}\\
0\ar[r] &\ker^S_{\Sigma}\ar[r]&\big(\Lambda[\Sigma_L]_0\big)_\Tht \ar[r]
       &\Lambda[\Sigma_L\smallsetminus S_L]_\Tht\,.
}
\]
Let $\iota\colon \Sigma_L\into \Sigma'_L$ denote the inclusion map.
Since by construction we have
\[\phi\circ\iota=\id\colon \Sigma_L\into\Sigma'_L\onto\Sigma_L,\]
we see that
\[\phi_*\circ \iota_*=\id\colon\hs \ker^S_\Sigma(\Lambda)\to\ker^S_{\Sigma'}(\Lambda)\to \ker^S_\Sigma(\Lambda).\]
Below we show that $\iota_*\circ\phi_*=\id$;
this will prove that $\iota_*$ is an isomorphism.

Let
\[
\lambdatil=\sum_{w\in \Sigma'_L} \lambda_w\cdot w \in \Lambda[\Sigma'_L]_0
\]
with $[\lambdatil]\in\ker^S_{\Sigma'}$.
Then $\sum_{w\in \Sigma_L'}\lambda_w=0$ and for each $v\in \Sigma'\smallsetminus S$ we have
\[\sum_{w\in \V_L(v)} \lambda_w\cdot w=\sum_{\vt\in\Theta} (\vt-1)\cdot \alpha_{v,\vt}\]
for some $\alpha_{v,\vt}\in \Lambda[\V_L(v)]$.
In particular, this holds for all $v\in\Sigma'\smallsetminus \Sigma$.

We have
\begin{align*}
(\iota_*\circ\phi_*)(\lambdatil)-\lambdatil
&=\sum_{v\in\Sigma'}\bigg( \phi\Big(\sum_{w\in\V_L(v)}\!\!\lambda_w\cdot w\Big)
   -\sum_{w\in\V_L(v)}\!\!\lambda_w\cdot w\bigg)\\
&=\sum_{v\in\Sigma'\smallsetminus\Sigma}\bigg( \phi\Big(\sum_{w\in\V_L(v)}\!\!\lambda_w\cdot w\Big)
   -\sum_{w\in\V_L(v)}\!\!\lambda_w\cdot w\bigg)
\end{align*}
where we may take the sum over $\Sigma'\smallsetminus\Sigma$
because for $v\in\Sigma$ and $w\in\V_L(v)$ we have $\phi(w)=w$.
It follows that
\begin{align*}
(\iota_*\circ\phi_*)(\lambdatil)-\lambdatil
&=\sum_{v\in\Sigma'\smallsetminus\Sigma}\phi\Big(\sum_{\vt\in\Theta}(\vt-1)\alpha_{v,\vt}\Big)
   -\sum_{v\in\Sigma'\smallsetminus\Sigma\ }\sum_{\vt\in\Theta}(\vt-1)\alpha_{v,\vt}\\
&=   \sum_{v\in\Sigma'\smallsetminus\Sigma}
    \ \sum_{\vt\in\Theta}(\vt-1)\cdot\big(\phi_*(\alpha_{v,\vt})-\alpha_{v,\vt}\big).
    \end{align*}
Write
\[\alpha_{v,\vt}=\sum_{w\in \V_L(v)}\!\! a_w\cdot w\in \Lambda[\V_L(v)]\subset \Lambda[\Sigma'_L] \]
and set
\[\beta_{v,\vartheta}=\phi_*(\alpha_{v,\vt})-\alpha_{v,\vt}=
         \sum_{w\in \V_L(v)}\!\! a_w\cdot \big(\phi(w)-w\big)\in \Lambda[\Sigma'_L]_0\hs.\]
Note that while $\alpha_{v,\vt}\in \Lambda[\Sigma'_L]$, we have
$\beta_{v,\vartheta}\in \Lambda[\Sigma'_L]_0$\hs.
Then
\[(\iota_*\circ\phi_*)(\lambdatil)-\lambdatil=
     \sum_{v\in\Sigma'\smallsetminus\Sigma}
     \ \sum_{\vt\in\Theta}(\vt-1)\cdot\beta_{v,\vt}
     =\ \sum_{\vt\in\Theta}\bigg((\vt-1)
     \cdot\!\!\!\!\sum_{v\in\Sigma'\smallsetminus\Sigma}\!\!\!\beta_{v,\vt}\bigg)
     \]
 with $\beta_{v,\vt}\in \Lambda[\Sigma'_L]_0$\hs, whence
\[ [(\iota_*\circ\phi_*)(\lambdatil)]-[\lambdatil]=0\in \big(\Lambda[\Sigma'_L]_0\big)_\Tht\]
and $(\iota_*\circ\phi_*)([\lambdatil])=[\lambdatil]$.
Thus
\[\iota_*\circ \phi_*=\id\colon\hs \ker^S_{\Sigma'}(\Lambda)\to\ker^S_{\Sigma}(\Lambda)\to \ker^S_{\Sigma'}(\Lambda)\]
and $\iota_*$ is an isomorphism, which completes the proof of the theorem.
\end{proof}

\begin{proposition}
Let $S_1\subseteq S_2\subset \V_K$ be finite subsets.
Let $\Sigma\subset \V_K$ be a sufficiently large finite subset
for $S_2$ (and hence for  $S_1$).
Then the following diagram with natural arrows commutes:
\[
\xymatrix{
\ker_\Sigma^{S_1}(\Lambda)\ar[r]^-\sim\ar[d]  & \Sha^{S_1}(\Lambda)\ar[d]\\
\ker_\Sigma^{S_2}(\Lambda)\ar[r]^-\sim  & \Sha^{S_2}(\Lambda)
}
\]
\end{proposition}

\begin{proof}
Both horizontal isomorphisms are induced by inclusions $\iota$, and the vertical arrows are the inclusions
$\ker_\Sigma^{S_1}(\Lambda)\subseteq \ker_\Sigma^{S_2}(\Lambda)$ and $\Sha^{S_1}(\Lambda)\subseteq  \Sha^{S_2}(\Lambda)$.
\end{proof}

\section{Algorithm}
\label{s:algorithm}

Let the finite group $\Theta=\Gal(L/K)$ be as above, and denote by
\[\Lambda=\pi_1^\alg(G)\]
the algebraic fundamental group of \cite{Borovoi-Memoir}.
We choose a finite system of generators  $e_1,\dots,e_\smm$ of $\Lambda$.
We have a homomorphism
\begin{equation*}
\vk\colon \Z^{\smm}\to\Lambda
\end{equation*}
given by $\vk(\varepsilon_i)=e_i$ for the standard basis $(\varepsilon_i)$ of $\Z^\smm$.
Then $\ker\vk$ is a  finitely generated subgroup of $\Z^{\smm}$;
let $\Nt$ be an $\tt nn\times\smm$ matrix whose rows generate $\ker\vk$.
We denote  $\Mt=I_\smm$ (the identity matrix); then $\Lambda= \Mta/\Nta$.

The finite group $\Theta$ acts on the finitely generated group $\Lambda$.
For each $\vt\in\Theta$ we choose an $\smm\times\smm$ integer matrix (not uniquely determined)
$\rho(\vt)=\big(\rho(\vt)_i^j\big)_{1\le i,j\le\smm}$ such that
\[\vt(e_i)=\sum_j\rho(\vt)_i^j\hs e_j\hs\]
where $\rho(\vt)_i^j$ is the entry in row $j$, column $i$.
In the program, we {\em define} $\Lambda$ by  integers $\smm$ and $\tt nn$
and an $\nn\times\smm$  matrix $\Nt$ as $\Lambda=\Mta/\Nta$.
We define an action of $\Theta$ on $\Lambda$ by a family $\big(\rho(\vt)\big)_{\vt\in\Theta}$
of $\smm\times\smm$ integer matrices.
Using our function {\tt IsAction}, we check that these data
indeed define a left action of $\Theta$ on $\Lambda$, that is,
\[\rho(\vt)\Nta\subseteq\Nta,\quad
(\rho(1)-I_\smm)\cdot \Z^\smm\subseteq\Nta,\quad\text{and}
\quad \big(\rho(\vt_1\vt_2)-\rho(\vt_1)\rho(\vt_2)\big)\cdot  \Z^\smm\subseteq\Nta\]
for all $\vt,\vt_1,\vt_2\in\Theta$, where we regard elements of $\Z^\smm$ as column vectors.

We assume that we are given a set $S$, and we choose a sufficiently large set
$\Sigma=S\sqcup (\Sigma\smallsetminus S)$ (Definition \ref{d:sl}).
In terms of places of $K$ this means:
\begin{equation}
\begin{aligned}\label{e:v-Sigma}
&\textit{for any decomposition subgroup $\Theta_v\subseteq \Theta$
for $v\in \V(K)\smallsetminus \Sigma$, there exists} \\
&\textit{$v_\Sigma\in \Sigma\smallsetminus S$
with $\Theta_{v_\Sigma}\supseteq \Theta_v$ up to conjugacy.}
\end{aligned}
\end{equation}
We assume that we know the decomposition subgroups $\Theta_v$ for all $v\in \Sigma$.
We construct a $\Theta$-set $\Sigma_L$ by
\[ \Sigma_L=\bigsqcup_{v\in\Sigma} \Theta/\Theta_v\hs.\]
Consider the group $\Lambda[\Sigma_L]$ on which the  finite group $\Theta=\Gal(L/K)$ and the group ring
$\Z[\Theta]$ act.
We wish to compute
\[ \ker_\Sigma^S(\Lambda) =
\ker\Big[ \big(\Lambda[\Sigma_L]_0\big)_\ThT\to \Lambda[\Sigma_L\smallsetminus S_L]_\Theta\Big].\]
We also consider the ``extended kernel''
\[ {\rm Eker}_\Sigma^S(\Lambda)\coloneqq
\ker\Big[ (\Lambda[\Sigma_L]_0)_\Theta\to \Lambda[\Sigma_L\smallsetminus S_L]_\Theta\Big];\]
then $\ker_\Sigma^S(\Lambda)=\big({\rm Eker}_\Sigma^S(\Lambda)\big)_{\rm Tors}\hs.$

By abuse of notation, we write $\Sigma$ for $\Sigma_L$ and $S$ for $S_L$.
We have
\[{\rm Eker}_\Sigma^S(\Lambda)=\ker \Bigg[\frac{\Lambda[\Sigma]_0}{I_\Theta\, \Lambda[\Sigma]_0}
\,\to\,
\frac{\Lambda[\Sigma\smallsetminus S]}{ I_\Theta\, \Lambda[\Sigma\smallsetminus S]}\Bigg] \]
where
\[ \Lambda[\Sigma]_0=\Big\{\sum \lambda_w\cdot w\in\Lambda[\Sigma]\ \big|\ \sum\lambda_w=0\Big\},\]
and $I_\Theta\subset \Z[\Theta]$
is the ideal generated by the elements $\vt-1$ for $\vt\in\Theta$.
We see that
\[{\rm Eker}_\Sigma^S(\Lambda)= \frac
{\Lambda[\Sigma]_0\cap \big(\Lambda[S]\oplus I_\Theta\, \Lambda[\Sigma\smallsetminus S]\big)}
{I_\Theta\,\Lambda[\Sigma]_0}.
\]

We pass to the preimages in $\Mta[\Sigma]$  of the three groups
\[\Lambda[\Sigma]_0,\quad I_\Theta\,\Lambda[\Sigma]_0,\quad
\Lambda[S]\oplus I_\Theta\,\Lambda[\Sigma\smallsetminus S] \]
under the canonical homomorphism $\Mta[\Sigma] \to \Lambda[\Sigma]$.
These preimages are:
\begin{gather*}
\Mta[\Sigma]_0+ \Nta[\Sigma],
\quad
I_\Theta\,\Mta[\Sigma]_0 + \Nta[\Sigma],
\quad
\Mta[S]\,\oplus\, I_\Theta\,\Mta[\Sigma\smallsetminus S] + \Nta[\Sigma].
\end{gather*}
We obtain that
\[ {\rm Eker}_\Sigma^S(\Lambda)\,\cong\,\frac
{\big(\,\Mta[\Sigma]_0 +\Nta[\Sigma]\,\big)\,\cap\,
\big(\,\Mta[S]\oplus I_\Theta\,\Mta[\Sigma\smallsetminus S]
+\Nta[\Sigma]\,\big)}
{ I_\Theta\,\Mta[\Sigma]_0+\Nta[\Sigma]}\,.
\]

It is easy to compute  three matrices
$\tt Z0, C0$, and $\tt CC$ in Magma such that
\begin{gather*}
\la{\tt Z0}\ra=\Mta[\Sigma]_0,\quad {\tt \la C0\ra}=  I_\Theta\,\Mta[\Sigma]_0,\quad
{\tt  \la CC\ra}= \Mta[S]\,\oplus\, I_\Theta\,\Mta[\Sigma\smallsetminus S].
\end{gather*}
Let $\tt NN$ denote the {\em diagonal join}
of $\#\Sigma_L$ copies of $\tt N$;
then $\Nta[\Sigma]={\tt\la NN\ra}$.
We compute three matrices
\begin{gather*}\tt Z0N:=VerticalJoin(Z0,NN),\quad
C0N:=VerticalJoin(C0,NN),\\
\tt CCN:=VerticalJoin(CC,NN).
\end{gather*}
Then
\begin{gather*}
{\tt \la Z0N\ra}=\Mta[\Sigma]_0 + \Nta[\Sigma],\quad
\la{\tt C0N}\ra=  I_\Theta\,\Mta[\Sigma]_0 + \Nta[\Sigma],\\
\la{\tt CCN}\ra= \Mta[S]\,\oplus\, I_\Theta\,\Mta[\Sigma\smallsetminus S]
+ \Nta[\Sigma].
\end{gather*}
Using our Magma function {\tt Meet}, we compute a matrix
\[\tt mZ0NCCN:=Meet(Z0N, CCN)\quad\
\text{such that}\quad\
\tt \la mZ0NCCN\ra=\la Z0N\ra\cap\la CCN\ra .
\]
Using our Magma function ${\tt ComplementIntMat}$ (written by Willem A. de Graaf),
we compute
\[\tt comp, bas, x := ComplementIntMat(\ mZ0NCCN,\  C0N\ ).\]
Here $\tt bas$ is a matrix with $n_{\tt bas}$ rows $\tt bas[i]$,
and $\tt x$ is a sequence of positive integers $\tt x[i]$ of length $n_{\tt bas}$
such that  ${\tt x[i]* bas[i]}, \ i=1,\dots,\ n_{\tt bas}\hs$,
constitute a basis of $\tt \la C0N\ra$
and that  the vectors $\tt bas[i]$ belong to $\tt \la mZ0NCCN\ra$
and represent a set of generators of the finite abelian group
\begin{equation*}
\tt \big(\la mZ0NCCN\ra/  \la C0N\ra\big)_\Tors\hs,
\end{equation*}
with each class $\tt [bas[i]]$ of order $\tt x[i]$.
Furthermore, we have an isomorphism of finite abelian groups
\[\tt \bigoplus_i \Z/x[i]\Z
\cong(\,{\tt\la mZ0NCCN\ra\, /\,  \la C0N \ra}\, )_{\rm Tors}\]
sending $[1]\in \tt \Z/x[i]\Z$ to the class $\tt [bas[i]\hs]$ of $\tt bas[i]$.
It follows that
\[\Sha^S(\Lambda)\cong \ker_\Sigma^S(\Lambda)
\cong (\,{\tt\la mZ0NCCN\ra\, /\,  \la C0N \ra}\, )_{\rm Tors}\cong \bigoplus_i \Z/x[i]\Z,\]
which computes $\Sha^S(\Lambda)$.

It remains to compute the homomorphisms
\begin{equation}\label{e:lvv} l_v\colon \Sha^S(\Lambda)\cong\ker_\Sigma^S(\Lambda)\to \Lambda_\Thvt
\end{equation}
for real places $v$ in $S$.
It suffices to compute the images
$l_v\tt\big[ bas[i]\big]$ of the generators
$\tt \big[ bas[i]\big]$.
We may regard each  $\tt bas[i]$ as an element of $\Mta[\Sigma]$,
say, ${\tt bas[i]}=\sum_{w\in\Sigma} m_w\cdot w$.

Let $v\in S_\infty\subset \Sigma$. Choose a place  $w_0$ of $L$ over $v$.
Recall that $\V_L(v)$ denotes the preimage of $v$ in $\V_L$.
For each $w\in \V_L(v)$,  choose $\vt_w\in\Theta$ such that $\vt_w\cdot w_0=w$,
and consider
\[ m_v\coloneqq\sum_{w\in \V_L(v)}\vt_w^{-1}\cdot m_w\hs.\]
We write $\Theta_0$ for the decomposition subgroup (stabilizer)
$\Theta_{w_0}$ of $w_0$ in $\Theta$.
Then $l_v\tt\big[bas[i]\hs\big]$
is the image in $\Lambda_{\Theta_0}$ of $m_v\in\Mta$
under the natural homomorphism
\[ \la M\ra=\Z^\smm\to \Lambda\to\Lambda_{\Theta_0}\hs.\]
Since $ \tt x[i]\cdot bas[i]\in \la C0N\ra$, we have
\[{\tt x[i]}\cdot [m_v]=0\in \Lambda_{\Theta_0}\hs,\]
whence $[m_v]$ is a torsion element, $[m_v]\in\Lambda_{\Theta_0\hs,\Tors}$\hs.
Thus we have described the homomorphism \eqref{e:lvv}.

\section{Examples: computation of the invariant factors of $\Sha^S(\Lambda)$}
\label{s:examples}

Let $F=\Q(\zeta_8)\subset \Qbar$,
where $\zeta_8$ is a primitive 8th root of unity.
Write $\Theta=\Gal(F/\Q)$.
We embed $\Z/8\Z$ into $(R_{F/\Q}\SL_8)(\Q)=\SL_8(F)$
by sending the generator $1+8\Z$ to the scalar matrix $\zeta_8\cdot I_8$.
We set $A=\Z/8\Z\subset R_{F/\Q}\SL_8$.
Following Sansuc \cite[Example 5.8]{Sansuc}, 
we consider the quasi-split, non-split $\Q$-group $G=(R_{F/\Q}\SL_8)/A$.
Write $\Lambda=\pi_1^\alg(G)$. By \cite[Corollary C.2.2]{BK}, we have
\[ \Lambda\cong\Hom(\X^*(A),\Q/\Z)\]
where $\X^*$ denotes the character group.
We obtain:
\begin{align*}
 &\X^*(A)=\Hom(\Z/8\Z,\Qbar^\times)\cong\mu_8\hs,\\
&\Lambda\cong\Hom(\mu_8,\Q/\Z)=\Hom(\mu_8,\tfrac18\Z/\Z)\cong\Hom(\mu_8,\Z/8\Z).
\end{align*}
We fix an isomorphism (not $\Theta$-equivariant)
 $\Hom(\mu_8,\Z/8\Z)\isoto\Z/8\Z$
by sending a homomorphism $\chi\colon\mu_8\to\Z/8\Z$ to $\chi(\zeta_8)\in \Z/8\Z$.
We identify $\Lambda$ with $\Z/8\Z=\{\bar x=x+8\Z\ |\  x\in \Z\}$.
Then $\Aut(\Lambda)=(\Z/8\Z)^\times=\{\hs\bar 1,\bar 3,\bar 5,\bar 7\hs\}$
acting on $\Lambda=\Z/8\Z$ by  $\bar a\cdot\bar x=\ov{ax}$.
The action of $\Theta=\Gal(F/\Q)$ on $\Lambda$ permits us to identify  $\Theta$ with $(\Z/8\Z)^\times$.
Namely, for $\bar a\in(\Z/8\Z)^\times$, 
the element $\vt_{\bar a}\in\Theta$ with $\vt_{\bar a}(\zeta_8)=\zeta_8^{\bar a}$
acts on $\Lambda=\Hom(\mu_8,\Z/8\Z)$ by $\chi\mapsto \chi\circ\vt_{\bar a}^{-1}$,
that is, by multiplication by $\bar a^{-1}=\bar a$.
The group $\Theta=(\Z/8\Z)^\times$ has the following subgroups:  four cyclic subgroups
$\langle\hs\bar1\hs\rangle,\langle\hs\bar3\hs\rangle, \langle\hs\bar5\hs\rangle, \langle\hs\bar7\hs\rangle$,
and one non-cyclic subgroup $\Theta$.

Let $K\subset \Qbar$ be a totally real field such that $K\cap F=\Q$,
and consider the compositum $L=K\cdot F=K(\zeta_8)$.
We identify the Galois groups $\Gal(L/K)\cong\Gal(F/\Q)=\Theta$.
We consider the group $G_K\coloneqq G\times_\Q K$.
Then $\pi_1^\alg(G_K)=\Lambda$ with $\Gal(\Qbar/K)$ acting on $\Lambda$
through the isomorphisms $\Gal(L/K)\cong\Gal(F/\Q)=\Theta$.

By \cite[Lemma 7.5.2]{BK}, for each  archimedean place $v$ of $K$,
the map $\ab\Lambda_v$ of \eqref{e:abLambda-v} is bijective.
(Indeed, for a real place $v$ of $K$, the $K_v$-group $G_{K_v}=G_K\times_K K_v\cong G\times_\Q\R$
is the same real group for every real place $v$ of every $K$ as above,
and \cite[Lemma 7.5.2]{BK} concerns only this real group.)
Therefore, by Corollary \ref{c:bijective}, for any finite subset $S\subset \V_K$
the map $\abLambda^S$ of \eqref{e:ab-Sha} is bijective as well.
Moreover, for each  archimedean place $v$ of $K$ and each place $w$ of $L$ over $v$,
we have $L_w/K_v\simeq\C/\R$. It follows that
the non-trivial element of $\Theta_w=\Gal(L_w/K_v)$ sends $\zeta$ to $\zeta^{-1}=\zeta^7$,
where $\zeta=\zeta_8$.
Hence $\Theta_v=\la\bar 7 \ra$.
In what follows we write $G$ for $G_K$.

\begin{example}\label{ex:cyclic}
We consider the case  $K=\Q\big(\sqrt{7}\hs\big)$.
By \cite[Lemma 7.5.1]{BK},
for this $K$  and the Galois extension $L=K(\zeta_8)$,
all decomposition groups are cyclic.
For each $a=1,3,5,7,$ we choose different finite places   $v_a$, $v'_a$, $v''_a$ of $K$
with decomposition group $\langle\hs\bar a\hs\rangle$
(this is possible by the Chebotarev density theorem).
Using our Magma program, we compute
the finite abelian group $\Sha^S(\Lambda)\cong\ker_\Sigma^S(\Lambda)$
with $S=\varnothing$, $S=\{v_a\}$ and $S=\{v_a,v'_b\}$ 
(with {\em different} $v_a$ and $v'_b$ even when $a=b$)
for each $a=1,3,5,7$ and $b=1,3,5,7$.
For our computation we take $\Sigma\smallsetminus S=\{v''_3,v''_5,v''_7\}$ 
with the places {\em different} from those appearing in $S$. 
Since all decomposition groups for $L/K$ are cyclic, 
condition \eqref{e:v-Sigma} is satisfied.
We compute the invariant factors of $\Sha^S(\Lambda)$,
which gives  the cardinality $\#\Sha^S(K,G)=\#\Sha^S(\Lambda)$
(because the  map $\ab\Lambda^S$ of \eqref{e:ab-Sha} is bijective).


\begin{longtable}{|r|l|l|l|l|l|}
\caption{$\Sha^S(\Lambda)\cong\Sha^S(K,G)$ for $K=\Q\big(\sqrt7\big)$}
\endfirsthead
\hline
\endhead
\hline
\endfoot
\endlastfoot

\hline
    &$\ \vn$    & \ $\la1\ra\hs$ &\ $\la3\ra\hs$ &\ $\la5\ra\hs$ &\ $\la7\ra\hs$
    \\
\hline
$\vn\ $   &   2    &   2,4   &   2   &   4   & 2
\\
\hline
$\ \la1\ra\ $
& 2,4 & 2,4,8 & 2,8 & 4,8 & 2,8
\\
\hline
$\ \la3\ra\ $
& 2 & 2,8 &  2,2 & 8 & 4
\\
\hline
$\ \la5\ra\ $
& 4 & 4,8 & 8 & 4,4  & 8
\\
\hline
$\ \la7\ra\ $
& 2 & 2,8  & 4  & 8 & 2,2
\\
\hline
\end{longtable}
In the cell with coordinates $\vn,\vn$ in this table,
the number $2$ is the invariant factor
of the finite abelian group $\Sha^S(\Lambda)$
for $S=\vn$.
In the cell with coordinates $\vn, \la1\ra$,
the numbers $2,4$ are the invariant factors for $S=\{v_1\}$.
In the cell with coordinates
$\la1\ra,\la1\ra$,
the numbers $2,4,8$ are the invariant factors
of the finite abelian group $\Sha^S(\Lambda)$
for $S=\{v_1,v_1'\}$ with {\em different} places
of $K$ having trivial decomposition groups
$\la\bar 1\ra,\la\bar 1\ra$.
In the cell with coordinates $\la3\ra,\la1\ra$,
the numbers $2,8$ are the invariant factors of the finite abelian group $\Sha^S(\Lambda)$
for  $S=\{v_3,v_1'\}$ with (different) places
of $K$ having  decomposition groups $\la\bar 3\ra,\la\bar 1\ra$,
and similarly for the other cells.

In particular, for $S=\vn$ we have $\Sha^\vn(K,G)=\Sha(K,G)\cong\Z/2\Z$,
which is well known; see \cite[Example 5.8]{Sansuc} and \cite[Formula (7.7) and Remark 7.5.4]{BK}.
\end{example}

We also consider the field $K=\Q(\alpha)$ where $\alpha$ is a root of the polynomial $f(x)=x^3+x^2-10x-8$.

\begin{proposition}\label{p:cubic}
For $K=\Q(\alpha)$ as above,
set $L=K(\zeta_8)$, $\Theta=\Gal(L/K)$. Then:
\begin{enumerate}
\item[\rm(i)] $K$ is totally real, with three real places.
\item[\rm (ii)] $\Theta\cong(\Z/8\Z)^\times$.
\item[\rm (iii)] $K$ has exactly three places $v$ with non-cyclic decomposition group,
namely, the three places over $2$ with decomposition group $\Theta$.
\end{enumerate}
\end{proposition}

\begin{proof}[Sketch of a proof]
The field $K$ is the cubic subfield of $\Q(\zeta_{31})$ (of discriminant $31^2$).
It is totally real, and we have $K\cap F=\Q$, since $[K:\Q]=3$ is odd.
Since  $2$ is a cube modulo $31$ ($2^{10}\equiv1$), we see that $2$ splits completely in $K$.
Hence, $K_v\cong\Q_2$ and $\Theta_v=\Gal(\Q_2(\zeta_8)/\Q_2)=\Theta$ for $v\,|\,2$.
The extension $L/K$ is unramified outside $2$ and $\infty$, and so all other decomposition groups are cyclic.
\end{proof}

\begin{example}\label{ex:Theta}
Let $K=\Q(\alpha)$ be as above.
For each $a=1,3,5,7,$ we choose two different finite places  $v_a$, $v'_a$ of $K$
with decomposition group $\langle\hs\bar a\hs\rangle$
(this is possible by the Chebotarev density theorem).
We choose three different places $v_\Theta, v'_\Theta, v''_\Theta$ of $K$
with decomposition group $\Theta$ 
(this is possible because by Proposition \ref{p:cubic}(iii)
the field $K$ has three places with non-cyclic decomposition group).
Using our Magma program, we compute
the finite abelian group $\Sha^S(\Lambda)\cong\ker_\Sigma^S(\Lambda)$
for each $S=S_\re\sqcup S_\other$, where $S_\re=\{v_7^{(1)},v_7^{(2)}\}$
consists of two different {\em real} places
(with decomposition subgroup $\la\bar 7\ra$), and where
$S_\other=\varnothing$, or $S_\other=\{v_a\}$, or  $S_\other=\{v_a,v'_b\}$
(with {\em different} $v_a$ and $v'_b$ even when $a=b$)
for each $a=1,3,5,7,\Theta$ and $b=1,3,5,7,\Theta$.
For our computation we take $\Sigma\smallsetminus S=\{v''_\Theta\}$
(with $v''_\Theta\notin S$); then condition \eqref{e:v-Sigma} is clearly satisfied.
We compute the invariant factors of $\Sha^S(\Lambda)$,
which gives the cardinality $\#\Sha^S(K,G)=\#\Sha^S(\Lambda)$
(because the map $\ab\Lambda^S$ of \eqref{e:ab-Sha} is bijective).


\begin{longtable}{|r|l|l|l|l|l|l|}
\caption{$\Sha^S(\Lambda)\cong\Sha^S(K,G)$ for $K=\Q(\alpha)$}
\endfirsthead
\hline
\endhead
\hline
\endfoot
\endlastfoot

\hline
    &$\ \vn$    & \ $\la1\ra\hs$ &\ $\la3\ra\hs$ &\ $\la5\ra\hs$ &\ $\la7\ra\hs$ & $\Theta$
    \\
\hline
$\vn\ $   &   2    &   2,8   &  2,2   &   2,4   & 2,2 & 2,2
\\
\hline
$\ \la 1 \ra\ $
& 2,8 & 2,8,8 & 2,2,8 & 2,4,8 & 2,2,8 & 2,2,8
\\
\hline
$\ \la 3 \ra\ $
& 2,2 & 2,2,8  & 2,2,2 & 2,2,4 & 2,2,2 & 2,2,2
\\
\hline
$\ \la 5 \ra\ $
& 2,4 & 2,4,8 & 2,2,4 & 2,4,4  & 2,2,4 & 2,2,4
\\
\hline
$\ \la 7 \ra\ $
& 2,2 & 2,2,8 & 2,2,2 & 2,2,4 & 2,2,2 & 2,2,2
\\
\hline
$ \Theta\ $
& 2,2  & 2,2,8 & 2,2,2 & 2,2,4 & 2,2,2 & 2,2,2
\\
\hline
\end{longtable}

The notation in Table 2 is similar to that in Table 1.
The row and column labels describe $S_\other$
where $S=S_\re\sqcup S_\other$.
For instance, the cell $(\vn,\vn)$ corresponds to $S=S_\re=\{v_7^{(1)},v_7^{(2)}\}$,
not to $S=\vn$.

\end{example}

\appendix

\section{A simple method of computing $\#\Sha^S(K,G)$}
\label{app:simple}

Here we describe a simple method of computing the cardinality of $\Sha^S(K,G)$
suggested by Zev Rosengarten (private communication).

\begin{proposition}
For a reductive group $G$ over a global field $K$,
there is a canonical exact sequence
\[1\to \Sha^\vn(K,G)\to H^1(K,G)\labelto{\loc}
\bigoplus_{v\in\V_K} H^1(K_v,G)\labelto\psi \Lambda_\Tht\]
where the map $\psi$ is the composite
\[\bigoplus_{v\in\V_K} H^1(K_v,G)\labeltoo{\bigoplus\ab\Lambda_v}
\bigoplus_{v\in\V_K}\Lambda_\Thvt\labelto\sigma\Lambda_\Tht\hs.\]
Here the maps $\ab\Lambda_v$ are as in Section \ref{s:reduction}, 
and $\sigma$ is the summation map:
\[\big(\hs[\lambda_v]_\Thv\big)_{v\in\V_K}\,\longmapsto
\,\sum_{v\in\V_K}\big[\lambda_v\big]_\Theta\]
where on the left-hand side, $[\lambda_v]_\Thv$ denotes the class in $\Lambda_\Thvt$\hs,
and on the right-hand side, 
$[\lambda_v]_\Theta$ denotes the class in $\Lambda_\Tht$\hs.
\end{proposition}

\begin{proof}
In the number field case, this is \cite[Theorem 5.15]{Borovoi-Memoir}
(it goes back to \cite[Corollary 2.5 and Proposition 2.6]{Kottwitz-86}\hs).
In the function field case, see \cite[Theorem 4.2]{Thang} or \cite[Theorem 1.6]{Rosengarten:Tamagawa}.
\end{proof}

\begin{corollary}\label{c:S-exact}
For a finite set $S$ of places of $K$, we have a canonical exact sequence
\[ 1\to\Sha^\vn(K,G)\to \Sha^S(K,G)\labelto{\loc_S}\prod_{v\in S} H^1(K_v,G)
  \labelto{\psi_S}\Lambda_\Tht\hs, \]
where the map $\psi_S$ is the composite
\[\psi_S\colon\, \prod_{v\in S} H^1(K_v,G) \labeltoo{\bigoplus\ab\Lambda_v}
\bigoplus_{v\in S}\Lambda_\Thvt\labelto{\sigma_S}\Lambda_\Tht\] 
and the summation map $\sigma_S$ is the restriction of the map $\sigma$.
\end{corollary}

\begin{proof}
It suffices to prove the exactness at $\prod_{v\in S}H^1(K_v,G)$. If $x\in\Sha^S(K,G)$, then by the proposition we have
$\psi_S(\loc_S x)=\psi(\loc\,x)=0$. Conversely, an element of
$\prod_{v\in S}H^1(K_v,G)$ contained in $\ker\psi_S$, when extended by trivial classes outside
$S$, lies in $\ker\psi$, and by the proposition it is $\loc(x)$ for some $x\in H^1(K,G)$. Since
$\loc_v(x)=1$ for $v\notin S$, we have $x\in\Sha^S(K,G)$, as desired.
\end{proof}

\begin{proposition}
If $S$ is as in Corollary \ref{c:bijective-Intro}, we have a canonical commutative diagram with exact rows and bijective vertical arrows
\[
\xymatrix@C=12mm{
 1\ar[r] &\Sha^\vn(K,G)\ar[r]\ar[d]^-{\ab^\vn} &\Sha^S(K,G)\ar[r]^-{\loc_S}\ar[d]^-{\ab^S} 
 &\prod_{v\in S} H^1(K_v,G)\ar[d]^-{\ab_S}  \ar[r]^-{\psi_S}&\Lambda_\Tht\ar@{=}[d]\\
1\ar[r] &\Sha^\vn_\ab(K,G)\ar[r] &\Sha^S_\ab(K,G)\ar[r]^-{\loc_S}&\bigoplus_{v\in S} H^1_\ab(K_v,G)
  \ar[r]^-{\psi_S}&\Lambda_\Tht
}
\]
where the bottom arrow $\psi_S$ is the homomorphism defined by
\[\psi_S\colon\, \bigoplus_{v\in S} H^1_\ab(K_v,G) \labeltoo{\bigoplus\vpL_v}
\bigoplus_{v\in S}\Lambda_\Thvt\labelto{\sigma_S}\Lambda_\Tht\hs.\]
\end{proposition}

\begin{proof}
The top row is exact by Corollary \ref{c:S-exact}.
The vertical arrow $\ab^\vn$ is bijective by Corollary \ref{c:no-infinite-Intro}.
Since $S$ is as in Corollary \ref{c:bijective-Intro},
the vertical arrow $\ab_S$ is bijective
(see the discussion after \eqref{e:abLambda-v}\hs),
and by that  corollary 
the vertical arrow $\ab^S$ is bijective too.
It follows that the bottom row is exact as well.
\end{proof}

From now on we assume that $S$ is as in Corollary \ref{c:bijective-Intro}.
By that corollary,
to compute $\#\Sha^S(K,G)$, 
it remains to compute the order (cardinality) of the abelian group $\Sha^S_\ab(K,G)$.
For any $v\in S$ we have an isomorphism of finite abelian groups $H_\ab^1(K_v,G)\isoto \Lambda_\Thvt$.
Thus the group $\Sha^S_\ab(K,G)$ fits into an exact sequence
\[0\to \Sha^\vn_\ab(K,G)\to \Sha^S_\ab(K,G)\to \bigoplus_{v\in S} \Lambda_\Thvt 
    \labelto{\sigma_S}\Lambda_\Tht\hs,\]
and the order of  $\Sha^S_\ab(K,G)$ is equal to the product of the orders of 
$\Sha^\vn_\ab(K,G)$ and of the kernel
\begin{equation}\label{e:ker}
\ker\Big[\bigoplus_{v\in S} \Lambda_\Thvt \labelto{\sigma_S}\Lambda_\Tht\Big].
\end{equation}
One can compute the group $\Sha^\vn_\ab(K,G)$, for instance, as in \cite[Subsection 4.12 and Corollary 5.13]{Borovoi-Memoir}:
\[\Sha^\vn_\ab(K,G)\cong\coker\Big[\bigoplus_{v\in \V_K} H_1(\Theta_v,\Lambda)\to H_1(\Theta,\Lambda)\Big]\]
where $H_1$ denotes group homology. One can compute the kernel \eqref{e:ker} from the definition.

\begin{example}\label{ex:Zev}
We use the notation of Example \ref{ex:cyclic}.
Let $G$ be the $\Q$-group of Section \ref{s:examples}.
Let $K=\Q\big(\sqrt7\big)$, and let $S=\{v_5\}$.
We write $G$ for $G_K$.
Then it is well known that $\Sha^\vn(K,G)\cong\Z/2\Z$; see, for instance,
Sansuc \cite[Example 5.8]{Sansuc}. 
Moreover, we have
\[\bigoplus_{v\in S} \Lambda_{\Theta_{v},\Tors}=\Lambda_{\Theta_{v_5},\Tors}\cong\Z/4\Z,\qquad\Lambda_\Tht\cong\Z/2\Z,\] 
and we have the homomorphism
\[\sigma_S\colon\,\Z/4\Z\cong \bigoplus_{v\in S} \Lambda_{\Theta_{v},\Tors}\, \longrightarrow\, 
\Lambda_\Tht\cong\Z/2\Z, \quad a+4\Z\mapsto a+2\Z \quad\text{for}\ a\in\Z.\]
We see that the kernel of $\sigma_S$ is of order $2$ and the order of $\Sha^\vn_\ab(K,G)$ is $2$,
whence the order of $\Sha^S_\ab(K,G)$ is $4$.
\end{example}

\begin{remark}
Note that the method of the body of the paper (based on \cite[Theorem 5.5.10]{BK}\hs) 
gives a bit more than the simple method of this appendix. 
Namely, the method of the appendix permits one to compute the isomorphism classes
of the subgroup  $\Sha^\vn_\ab(K,G)$ 
and the corresponding quotient group $\Sha^S_\ab(K,G)/\Sha^\vn_\ab(K,G)$, 
whereas the method of the body of the paper 
permits one to compute the isomorphism class of the group $\Sha^S_\ab(K,G)$ itself.
For instance, in Example \ref{ex:Zev} we see that $\Sha^S_\ab(K,G)$
fits into an exact sequence
\[0\to\Z/2\Z\to \Sha^S_\ab(K,G)\to\Z/2\Z\to 0,\]
whereas in Example \ref{ex:cyclic} we obtained in the cell $(\hs\la5\ra,\vn\hs)$
of Table 1 that  $\Sha^S_\ab(K,G)\simeq\Z/4\Z$.
\end{remark}


\section{Magma program}
\label{app:Magma}

Our program computes the group $\ker_{\Sigma}^S(\Lambda)\cong \Sha^S(\Lambda)$ of Theorem \ref{t:Sigma'-Intro}
with $\Lambda=\la M\ra/\la N\ra $.
Here $\Theta$ is a finite group of order $\tt gg$;
it is denoted by $G$ in the program.
Moreover, $\tt\la M\ra= \Z^\smm$ is a free
finitely generated abelian group with $\smm$ generators,
and $\tt\la N\ra$ is a subgroup of $\Z^\smm$ with $\tt nn$ generators
that are the rows of an $\tt nn\times \smm$ matrix $\tt N$.

Furthermore, $\Sigma= S\sqcup (\Sigma\smallsetminus S)$
is represented in the program by the finite sequence $(\Thv)_{v\in\Sigma}$ of {\em subgroups of $\Theta$}.
The result depends only on the conjugacy classes of the subgroups.
Actually, the result depends only on the multiplicity of each conjugacy class.
Moreover, the result does not depend
on the positive multiplicities for $\Sigma\smallsetminus S$,
and so we take each conjugacy class of subgroups at most once in $\Sigma\smallsetminus S$.
Returning to our field extension $L/K$,
we assume that \eqref{e:v-Sigma} holds.

In addition to $\ker_{\Sigma}^S(\Lambda)$, our program computes also
the homomorphisms
\[l_v\colon \ker_{\Sigma}^S(\Lambda) \to \Lambda_{\Theta_v,\Tors}\]
for $v\in S_\re$; see \eqref{e:lvv}.
Here $S_\re$ is the subsequence of $S$ corresponding to the real places of $K$;
then $\#\Theta_v=2$ for all $v\in S_\re$.

For explicit calculations we choose an example with $\smm=1$.
We take $\tt \la M\ra=\Z$, $\tt nn=1$, $\tt N=[[8]]$;
then $\tt\la N\ra= 8\Z$ and  $\Lambda=\Z/8\Z$.
We take $\Theta=\Aut \Lambda=(\Z/8\Z)^\times=\{\bar 1,\bar3,\bar5,\bar 7\}$
where for $a\in\Z$ we denote  $\bar a=a+8\Z$, the image of $a$ in $\Z/8\Z$.

\begin{itemize}
\item The elements $1,\, 2,\,3,\,4$ of $G=\{1,\dots,{\tt gg}\}$ correspond to 
$\bar 1,\bar3,\bar5,\bar 7$.

\item The subgroups $\tt HH[1], \dots, HH[5]$ correspond to 
$\la\bar 1\ra,\la\bar3\ra,\la\bar5\ra,\la\bar7\ra,\Theta$.

\item $\tt k\_1re,\ k\_1other$, and $\tt k\_2$ are sequences of indices in $\tt HH$
describing $S_\re$, $S\smallsetminus S_\re$, and $\Sigma\smallsetminus S$.

\item The output $\tt nOne$ is the list of invariant factors $\ne 1$ 
of $\ker_\Sigma^S(\Lambda)$,
and $\tt ellV$ gives the values of $l_v$ on the generators.
\end{itemize}

\bigskip

\centerline{\sc The program}
\medskip

\begin{lstlisting}[language=C++]

/* This program computes Sha^S(G,M/N) where G is a finite group of order gg acting on a finitely generated abelian group M/N, where M is the free abelian group  Z^mm with mm generators, and N is a subgroup of M with nn generators.
*/

ComplementIntMat:= function( full, sub )
/* This is a Magma version, written by Willem A. de Graaf, of the GAP function with the same name.
*/
     full:= Matrix( Integers(), full );
     H:= HermiteForm( full );
     full:= [ ];
     for i in [1..NumberOfRows(H)] do
         if not IsZero( H[i] ) then
            Append( ~full, Eltseq( H[i] ) );
         end if;
     end for;
     full:= Matrix( Integers(), full );

     n:= NumberOfColumns(full);
     m:= NumberOfRows( full );
     Vn:=KModule( Rationals(), n );
     W:= KModuleWithBasis( [ Vn!Eltseq(full[i]) : i in [1..m] ] );
     cfs:= [ Coordinates( W, Vn!Eltseq(sub[i]) ) : i in [1..NumberOfRows(sub)] ];
     A:= Matrix( Integers(), cfs );

     S,P,Q:= SmithForm( A );
     Qi:= Q^-1;
     cfs:= [ Eltseq(Qi[j]) : j in [1..m] ];
     vecs:= [ &+[ c[i]*(Vn!Eltseq(full[i])) :  i in [1..m] ] : c in cfs ];
     x:= [ ];
     minRC:=Minimum(NumberOfRows(S),NumberOfColumns(S));
     for i in [1..minRC] do
         if S[i][i] ne 0 then
            Append( ~x, S[i][i] );
         end if;
     end for;

     comp:= [ ];
     for i in [#x+1..#vecs] do
        Append( ~comp, Eltseq(vecs[i]) );
     end for;

     return Matrix( Integers(), comp),
	 Matrix( Integers(), [ Eltseq(vecs[i]) : i in [1..#x] ] ), Vector( Integers(), x );

end function;



IsGroup := function(gg, mu)
/* This function checks that the pair (G,mu), where G={1,..,gg} and mu: G x G -> G, is a group, and if yes, computes the inversion map iota.
*/
	boo:=true;	
	
	for g:=1 to gg do
		boo := boo and (mu[g,1] eq g) and (mu[1,g] eq g);
//      the unit condition
	end for;	
				
	for g1:=1 to gg do
	for g2:=1 to gg do
		for g3:=1 to gg do
			p12:=mu[g1,g2]; p23:=mu[g2,g3];
			boo:=boo and (mu[p12,g3] eq mu[g1,p23]);
//			 the associativity condition
		end for;
	end for;
	end for;
		
	iota:=[];

	for g:=1 to gg do
		boog:=false;
		for h:=1 to gg do
			boo1:=(mu[g][h] eq 1);
			if boo1 then
		    	iota[g]:=h;
				boog:=true;
				break;
			end if;
		end for;
		boo:=boo and boog;
//		 existence of inverse
	end for;
		
	return boo,iota;
	
end function;



IsAction:= function(gg,mu,mm,rho,nn,N)
/* This function, for a matrix N with mm columns and a family \rho(g) of mm x mm matrices, checks whether rho induces a left action of (G,mu) on Z^mm / <N> .
*/	
	boo:=true;	
	Zv:=RModule(Integers(), mm);
	vv:= Zero(Zv);			

	for g:=1 to gg do
	for i:=1 to nn do
		for j:=1 to mm do
			vv[j]:=N[i][j];
		end for;
		vv:=vv * Transpose(rho[g]);	
		ic:=IsConsistent(N, vv);
//		rho preserves N		
		boo:=boo and ic;
	end for;
	end for;
	
	Idm := ScalarMatrix(Integers(), mm, 1 );
	Dif:= Transpose (rho[1]-Idm);
	for i:=1 to mm do
			for j:=1 to mm do
				vv[j]:=Dif[i][j];
			end for;	
			ic:=IsConsistent(N, vv);
// 			rho[1] is 1  modulo N
			boo:=boo and ic;
		end for;
	
	for g:=1 to gg do
	for h:= 1 to gg do
		Dif:=Transpose( rho[mu[g][h]] - rho[g]*rho[h] );
		for i:=1 to mm do
			for j:=1 to mm do
				vv[j]:=Dif[i][j];
			end for;	
			ic:=IsConsistent(N, vv);
// 			rho is a well-defined left action modulo N
			boo:=boo and ic;
		end for;
	end for;
	end for;
		
	return boo;
	
end function;
			


Meet := function(A,B) //
// This function computes a matrix R such that <R> = <A> \cap <B>.

	AB:=VerticalJoin(A,B);
	NS:=NullspaceMatrix(AB);
	rA:=NumberOfRows(A);
	rNS:=NumberOfRows(NS);
	R:=ExtractBlockRange(NS, 1, 1, rNS,rA) * A;
	
	return R;
	
end function;
	
	
	
Mult:= function(g,X,mu)
// Left action of an element g of G on a subset X of G.	
	return { mu[g][x] : x in X };
end function;



CheckSubGroup:= function(G,H,mu,iota)
	
	boo:=true;
	for h1 in H do
		boo1:=iota[h1] in H;
	if not boo1 then
		boo:=false; break;
	else
		for h2 in H do
			boo2:= mu[h1][h2] in H;
			if not boo2 then
				boo:=false; break;
			end if;
		end for;
	end if;
	end for;
	
	return boo;
	
end function;



Cosets:=function(gg,H,mu)
/* This function computes the left cosets of G={1..gg} by a subgroup H.
*/
	cosets:=[]; ic:=0;
	booG:=[true : g in [1..gg] ];
	
for g:=1 to gg do
	if booG[g] then
		ic:=ic+1;
		cosets[ic]:= Mult(g,H,mu);
		for x in cosets[ic] do
			booG[x]:=false;
		end for;
	end if;
end for;

return cosets;

end function;



Num:= function(X, cosets, yy)
// The number (sic!) of a given left coset X of G.

	for i:= 1 to yy do
		if X eq cosets[i] then
			numX:=i; break;
		end if;		
	end for;
	
	return numX;
	
end function;
	
	
	
NotOne:= function(bas,x)

/* This function computes the sequence nOne of the elements of the integer vector x different from 1, and the matrix of the corresponding rows of the matrix bas
*/
	lbas:=NumberOfColumns(bas);
	d:=#Eltseq(x);
	gen:=[];
	s:=[];
	l:=0;
	
	for i:=1 to d do
		if x[i] ne 1 then
			Append(~s,x[i]);
			l := l+1;
			Append(~gen,bas[i]);
		end if;
	end for;
	
	gen := Matrix(Integers(), l, lbas,  gen);
	
	return  s, gen;
end function;



// ######################################

print "\n\nINPUT";	

//INPUT: G (that is, mu),  two integers gg, mm.
	
gg:=4;
	
//  A list of lists mu[1..gg][1..gg]  (the multiplication law)
	 mu:=[[1,2,3,4],[2,1,4,3],[3,4,1,2],[4,3,2,1]];		

// INPUT of Lambda, that is, mm, N, and rho
	 	
mm:=1;	


        
// #######################

rho:= [ ZeroMatrix(Integers(),mm,mm) :
    u in [1..gg] ];

rho[1]:=[[1]];
rho[2]:=[[3]];
rho[3]:=[[5]];
rho[4]:=[[7]];

//	INPUT: N
	
nn:=1;

N:=ZeroMatrix(Integers(),nn, mm);
N[1][1]:=8;


print "\ngg=", gg,
    "\nthe multiplication matrix: mu\n", mu;	
	
print "\nmm=", mm, "<M> = <I_mm>= Z^mm";

print "\nthe action of G: rho \n", rho;	

print "\nnn=", nn,
    "\nthe relations in <M> = Z^mm:  N", N;

print "\nOur G-module is Z^mm / <N>";	
	
boo, iota := IsGroup(gg,mu);	
print "\nIs G a group?  ", boo;
print "\nthe inversion map: iota ", iota;
		
print "\nIs rho an action of G on Z^mm / <N> ?  ",
    IsAction(gg,mu,mm,rho,nn, N);



// ##################################

print "\n\n\nCOMPUTING Sha^S(K,G)";



// G = {1... gg}

// The subgroups of G: HH[i]
HH:=[];

HH[1]:={1};
HH[2]:={1,2};
HH[3]:={1,3};
HH[4]:={1,4};
HH[5]:={1,2,3,4};

mHH := [];
mHH[1]:=[1];
mHH[2]:=[1,3];
mHH[3]:=[1,5];
mHH[4]:=[1,7];
mHH[5]:=[1,3,5,7];

boo := true;
for HHH in HH do
	print "HHH = ", HHH;
	boo := boo and  CheckSubGroup( gg, HHH, mu, iota );
end for;
print "All subsets HH[i] are subgroups: ", boo;



ShaS:= function(gg, mu, iota, mm, N, rho, HH,
    k_1re, k_1other, k_2 )

/* This is the main function of the program. It computes Sha^S=ker_\Sigma^S. Here k_1re is the sequence of real places in S, k_1other is the sequence of other places in S, and k_2 is the set of places in  Sigma - S.
*/



nn:= NumberOfRows(N);

k_1 := k_1re cat k_1other;

kk_1:=#k_1;
kk_2:=#k_2;


k:=k_1 cat k_2; // Sigma
kk:=#k;  // The number of places in Sigma


xx:=[ 0 : u in [1..kk]];

for i:=1 to kk do
	xx[i]:=gg div #HH[k[i]];
    // The number of places of L
    // over the place v[i] of K.
end for;



// #########################

//Computing the matrices A[g][i] of action of g on the left cosets modulo HH[k[i]]

AA:=[ [* ZeroMatrix(Integers(), xx[i], xx[i]) :
    i in [1..kk] *] : g in [1..gg] ];

for i:=1 to kk do
for g:=1 to gg do
	AA[g][i]:=ZeroMatrix(Integers(), xx[i], xx[i]);
	cosets:=Cosets(gg,HH[k[i]],mu);
	for x:=1 to xx[i] do
//		print g, i, x, xx[i];
		y:=Num(Mult(g,cosets[x],mu), cosets, xx[i]);
		AA[g][i][x][y]:=1;
	end for;
end for;
end for;


/* Here AA[g][i] is the permutation matrix of the action of g on the places of L lying over the place v[i] of K.
*/

    
// ###############################

// Computing the matrices Z0 and C0 for i:=1 to kk for zero-sum vectors

Z0:=[];
zz:=0;

C0:=[];
cc:=0;

lx:=0;
dd:=0;

for i:=1 to kk do
for x:= 1 to xx[i] do
	lx:=lx+1;	
for m:=1 to mm do
	dd:=dd+1;	
		
	B:=[* ZeroMatrix(Integers(), xx[i], mm) : i in [1..kk] *];
	Bp:=[* ZeroMatrix(Integers(), xx[i], mm) : i in [1..kk] *];

	B[i][x][m]:=1;

	if (i eq 1) and (x eq 1) then
		B[1][1][m]:=0;
	else
		B[1][1][m]:=-1;
	end if;
	
	z0:=&cat[ElementToSequence(B[j]) : j in [1..kk] ];
	Append(~Z0,z0); // The matrix of generators of the group
    // of zero-sum vectors in k_1 and k_2
	zz:=zz+1;
	
	for g:=1 to gg do
		Bp:=B;		
		Bp[i]:=B[i]-Transpose(AA[g][i])*B[i]*Transpose(rho[g]);
		Bp[1]:=B[1]-Transpose(AA[g][1])*B[1]*Transpose(rho[g]);
		b:=&cat[ElementToSequence(Bp[j]) : j in [1..kk] ];
		Append(~C0,b);	// The matrix of x - g.x
        // for zero-sum vectors x	in k_1 and k_2
		cc:=cc+1;
	end for;
	
end for;
end for;
end for;
	
	
Z0Mat:=Matrix(Integers(), zz, dd, Z0 );
C0Mat:=Matrix(Integers(), cc, dd, C0 );
// Two matrices with dd = mm * kk columns


for l:=1 to lx do
	if l eq 1 then
		NN:=N;
	else
		NN:=DiagonalJoin(NN,N);
	end if;
end for;
// NN is a matrix with dd columns: <NN> = <N>[Sigma]

Z0N:=VerticalJoin(Z0Mat,NN); // <Z0N> = Z0Mat> + <NN>
C0N:=VerticalJoin(C0Mat,NN); // <C0N> = C0Mat> + <NN>



// ################################################

// Computing the matrix CC for S for i:=1 to kk_1

CC:=[];

lx:=0; dd:=0; cc:=0;

for i:=1 to kk_1 do
for x:= 1 to xx[i] do
	lx:=lx+1;	
	for m:=1 to mm do
		dd:=dd+1;	
		
		B:=[* ZeroMatrix(Integers(), xx[i], mm) :
            i in [1..kk] *];
		Bp:=[* ZeroMatrix(Integers(), xx[i], mm) :
            i in [1..kk] *];

		Bp[i][x][m]:=1;
		b:=&cat[ElementToSequence(Bp[j]) : j in [1..kk] ];
		Append(~CC,b);
            //The matrix of ALL vectors in S (k_1).	
		cc:=cc+1;
	end for;
end for;
end for;


// ############################################

// Computing CC for Sigma - S for i:= kk_1+1 to kk_1+kk_2

for i:=kk_1+1 to kk_1+kk_2 do
for x:= 1 to xx[i] do
	lx:=lx+1;	
for m:=1 to mm do
	dd:=dd+1;	
		
	B:=[* ZeroMatrix(Integers(), xx[i], mm) :
        i in [1..kk] *];
	Bp:=[* ZeroMatrix(Integers(), xx[i], mm) :
        i in [1..kk] *];

	B[i][x][m]:=1;
	
	for g:=1 to gg do
		Bp:=B;
		Bp[i]:=B[i]-Transpose(AA[g][i])*B[i]*Transpose(rho[g]);
		b:=&cat[ElementToSequence(Bp[j]) : j in [1..kk] ];
		Append(~CC,b);	// The matrix of x - g.x
             // for all vectors x in k_2.
		cc:=cc+1;
	end for;
end for;
end for;
end for;


CCMat:=Matrix(Integers(), cc, dd, CC );	
CCN:=VerticalJoin(CCMat,NN);
/*
CCN is the set of all cocycles (=cochains)  in M[Sigma] (modulo <N>)
( where Sigma = k_1 \cup k_2 )
that become coboundaries over Sigma - S = k_2
*/


// ####################################

mZ0NCCN := Meet(Z0N,CCN);
// <mZ0NCCN> = <Z0N> \cap <CCN>
comp, bas, x:=ComplementIntMat( mZ0NCCN, C0N );
/* x =
=  ( ( <Z_0N> \cap <CCN> ) / <C0N> )_Tors
= ker_\Sigma^S
*/



if k_2 eq [2,3,4] then
sk_2:= "CYCLIC: [ [1,3], [1,5], [1,7] ]";
else sk_2:="THETA: [ [ 1, 3, 5, 7 ] ]";
end if;

// Printing the list of the decomposition
// subgroups in k_1re,  k_1other, and k_2.

Pk_1re:= "[ ";
for i in k_1re do
	Pk_1re:= Pk_1re cat Sprint(mHH[i]) cat "  ";
end for;
Pk_1re:=Pk_1re cat "]";

Pk_1other:= "[ ";
for i in k_1other do
	Pk_1other:= Pk_1other cat Sprint(mHH[i]) cat "  ";
end for;
Pk_1other:=Pk_1other cat "]";


print "\nk_1re:", Pk_1re,  "\nk_1other:", Pk_1other,
"\nS = k_1 := k_1re cat k_1other";
// [ mHH[i] : i in k_1 ];
print "Sig - S = k_2:", sk_2, "\n";
// [ mHH[i] : i in k_2 ];

nOne, gen:= NotOne(bas,x);

return nOne, gen, xx;

end function;

//###########################################



ellV := function(gg, mu, iota, mm, N, rho, HH, k_1re, l, gen, xx)

/* This function computes the homomorphisms l[v] for the real places in v in S_\infty = k_1re
*/
	if l eq 0 then return []; end if;
	
	mm_ellV := [];
	kk_1re := #k_1re;
	lv:=0;
	lw:=0;
	
	for v := 1 to kk_1re do
		lv:=lv+1;
		matV := [];
		for i := 1 to l do
			cosets := Cosets(gg, HH[k_1re[v]], mu );
			lw:=0;
			for vv:= 1 to v-1 do
				lw := lw+xx[vv];
			end for;
			
			for y := 1 to xx[v] do
				lw := lw+1;
				vecy:=ZeroMatrix(Integers(), 1, mm);
				for g := 1 to gg do
					yNum := Num( Mult( g, HH[k_1re[v]], mu ), cosets, xx[v] );
					if yNum eq y then
						gy := g; break;
					end if;
				end for ; // g	
				
			vec:=Matrix(Integers(),1, mm,
			[ gen[i][ (lw-1)*mm+j ] : j in [1..mm] ]);
			vecy := vecy+ vec*Transpose(rho[iota[gy]]);
			end for;  // y
		
			Append( ~matV, Eltseq(vecy) );
		end for; // i
		Append( ~mm_ellV, matV);
	end for; // v
	R_ellV:=mm_ellV;
	
	if mm eq 1 then
		R1 := ZeroMatrix( Integers(), l, kk_1re );
		
		for v := 1 to kk_1re do
			for i := 1 to l do
				R1[i][v] :=mm_ellV[v][i][1];
			end for;
		end for;
		R_ellV:=R1;
	end if;
		
	return R_ellV;
	
end function;
				
	

Running := procedure(k_1re, k_1other, k_2)

/* This procedure runs the main functions ShaS and ellV and prints the results.   
Here Lambda_{Theta_v,Tors}=Z/2Z for Theta_v=<7>, and so we reduce ellvE[i][j] modulo 2.
*/
	nOne, gen, xx := ShaS(gg, mu, iota, mm, N, rho, HH,
        k_1re, k_1other, k_2 );
	print "nOne=", nOne;	
	l := #nOne;
	print "gen = ",  gen;

    if k_1re ne [] then
        if l gt 0 then
            ellvE := ellV( gg, mu, iota, mm, N, rho, HH, k_1re, l, gen, xx );
            for i:= 1 to NumberOfRows(ellvE) do
                for j:= 1 to NumberOfColumns(ellvE) do
                    ellvE[i][j] := ellvE[i][j] mod 2;
                end for;
            end for;
            print "ellV =", ellvE;
        else
            print "ellV = [ ]   (ker_Sigma^S(Lambda) is trivial)";
        end if;
    end if;
	print "*********************";
    
end procedure;



//##########################
// THE MAIN CYCLES

// (1) : Below we compute Sha^S(Lambda) for K = Q(\sqrt{7}) and $L=K(\zeta_8)
// as in Example 5.1.

print "\nCase of  K = Q(sqrt{7}) \n\n";


k_1re := [];  // The set of *real* places in S
kk_1re := 0;


k_2:= [2,3,4];
sk_2:= "CYCLIC";
	
		

	print "########################################\n";
	print "## k_1 is EMPTY,  k_2 is ", sk_2, "  ##";
	
	k_1other:= [];
	Running(k_1re, k_1other, k_2);
	

	print "#########################################\n";
	print "## k_1 contains ONE place,  k_2 is ", sk_2, "  ##";
	
	for ii := 1 to 4 do
		
		k_1other := [ii];
		Running(k_1re, k_1other, k_2);
		
	end for;
	

	print "#############################################\n";
	print "## k_1 contains TWO places,  k_2 is ", sk_2, "  ##";

	for ii := 1 to 4 do
	for jj := 1 to 4 do
	
		k_1other := [ii,jj];
		Running(k_1re, k_1other, k_2);
		
	end for;
	end for;


/* (2) : Below we compute Sha^S(Lambda) for K = Q(alpha) and $L=K(zeta_8)
 as in Example 5.3. Here alpha  is a root of the polynomial
f(x) = x^3 + x^2 - 10*x - 8.  */

print "\nCase of  K = Q(alpha)  where  alpha  is a root \nof the polynomial",
"f(x) = x^3 + x^2 - 10*x - 8.\n\n";

k_1re := [ 4, 4 ];  //  The sequence of real places in S
kk_1re := #k_1re;

print "k_1re contains two places: ", k_1re;

booRe:= true;
for v in k_1re do
//	print "v = ", v, "   HH[v] = ", HH[v];
  	if #HH[v] ne 2  then
  		print "HH[v] is not real for v = ", v;
  		booRe := false;
  		break;
  	end if;
end for;
print "All subgroups  HH[v] for v in k_1re are of order 2: ", booRe;

k_2 := [ 5 ];
sk_2:="THETA";
		

	print "########################################\n";
	print "## k_1other is EMPTY,  k_2 is ", sk_2, "  ##";
	
	k_1other:= [];
	Running(k_1re, k_1other, k_2);
	

	print "#########################################\n";
	print "## k_1other contains ONE place,  k_2 is  ", sk_2, "  ##";
	
	for ii := 1 to 5 do
		
		k_1other := [ii];
		Running(k_1re, k_1other, k_2);
		
	end for;
	

	print "#############################################\n";
	print "## k_1other contains TWO places,  k_2 is ", sk_2, "  ##";

	for ii := 1 to 5 do
	for jj := 1 to 5 do
	
		k_1other := [ii,jj];
		Running(k_1re, k_1other, k_2);
		
	end for;
	end for;


\end{lstlisting}


\newcommand{\sk}{\smallskip}


\begin{thebibliography}{BdG24}

\bibitem[BP90]
{BP}
Armand Borel and Gopal Prasad,
{\em Addendum to: `Finiteness theorems for discrete subgroups of bounded covolume in semi-simple groups'}, 
Inst. Hautes \'Etudes Sci. Publ. Math. No. {\bf 71} (1990), 173--177.
\sk

\bibitem[BS64]
{BS}
Armand Borel and Jean-Pierre Serre,
{\em Th\'eor\`emes de finitude en cohomologie galoisienne},
Comment. Math. Helv. {\bf 39} (1964), 111--164.
\sk 


\bibitem[Bor98]
{Borovoi-Memoir}
Mikhail Borovoi,
{\em Abelian Galois cohomology of reductive groups},
Mem. Amer. Math. Soc. {\bf 132} (1998), no. 626.
\sk

\bibitem[Bor22]
{Borovoi-real}
Mikhail Borovoi, 
{\em Galois cohomology of reductive algebraic groups over the field of real numbers},
Commun. Math. {\bf 30} (2022), no. 3, 191--201.
\sk

\bibitem[BdG24]
{BdG}
Mikhail Borovoi and Willem A. de Graaf,
{\em Computing Galois cohomology of a real linear algebraic group},
J. Lond. Math. Soc. (2)  {\bf 109} (2024), no. 5, Article ID e12906, 53 pp.
\sk

\bibitem[BK26]
{BK}
Mikhail Borovoi and Tasho Kaletha,
with an appendix by Vladimir Hinich,
{\em Galois cohomology of reductive groups over global fields,}
Mem. Amer. Math. Soc. {\bf 319} (2026), no. 1621, vii+110 pp.
\sk

\bibitem[BT24]
{BT}
Mikhail Borovoi and Dmitry A. Timashev,
{\em Galois cohomology and component group of a real reductive group},
Israel J. Math. {\bf 259} (2024), no. 2, 937--1000.
\sk

\bibitem[CT08]
{CT}
Jean-Louis Colliot-Th\'el\`ene,
{\em R\'esolutions flasques des groupes lin\'eaires connexes},
J. Reine Angew. Math. {\bf 618} (2008), 77--133.
\sk

\bibitem[Ko84]
{Kottwitz-84}
Robert E. Kottwitz,
{\em Stable trace formula: cuspidal tempered terms},
Duke Math. J. {\bf 51} (1984), no. 3, 611--650.
\sk

\bibitem[Ko86]
{Kottwitz-86}
Robert E. Kottwitz,
{\em Stable trace formula: elliptic singular terms},
Math. Ann. {\bf 275} (1986), no. 3, 365--399.
\sk

\bibitem[Nis79]
{Nisnevich}
Yevsey A. Nisnevich, 
{\em A properness theorem for Galois cohomology over function fields} (in Russian),
Dokl. Akad. Nauk BSSR {\bf 23} (1979), no. 12, 1065--1068, 1148.
\sk

\bibitem[Ros21]
{Rosengarten:Tamagawa}
Zev Rosengarten, 
{\em Tamagawa numbers and other invariants of pseudoreductive groups over global function fields},
Algebra Number Theory {\bf 15} (2021), no. 8, 1865--1920.
\sk

\bibitem[Ros23]
{Rosengarten:Tate}
Zev Rosengarten, 
{\em Tate duality in positive dimension over function fields},
Mem. Amer. Math. Soc. {\bf 290} (2023), no. 1444, v+217 pp.
\sk


\bibitem[Rot09]
{Rotman}
 Joseph J. Rotman,
{\em An introduction to homological algebra,
second edition},
Universitext,
Springer, New York, 2009.
\sk

\bibitem[San81]
{Sansuc}
Jean-Jacques Sansuc, \emph{Groupe de {B}rauer et arithm\'etique des groupes
alg\'ebriques lin\'eaires sur un corps de nombres}, J. Reine Angew. Math.
{\bf 327} (1981), 12--80.
\sk

\bibitem[Tat66]
{Tate}
John Tate,
{\em The cohomology groups of tori in finite Galois extensions of number fields},
Nagoya Math. J. {\bf 27} (1966), 709--719.
\sk

\bibitem[Tha19]
{Thang}
Nguy\^{e}\~{n} Qu\^{o}\'{c} Th\v{a}\'{n}g, 
{\em On Galois cohomology of connected reductive groups and Kottwitz exact sequence},
Bull. Sci. Math. {\bf 151} (2019), 66--138.

\end{thebibliography}
\end{document}